\documentclass{amsart}
\usepackage{latexsym}
\usepackage{graphicx}
\usepackage[centertags]{amsmath}
\usepackage{amsfonts}
\usepackage{amssymb}
\usepackage{amsthm}
\usepackage{newlfont}
\usepackage{subfigure}
\usepackage{epsfig}

\theoremstyle{plain}
\newtheorem{theorem}{Theorem}[section]

\newtheorem{corollary}[theorem]{Corollary}
\newtheorem{lemma}[theorem]{Lemma}

\theoremstyle{definition}
\newtheorem{defn}[theorem]{Definition}
\newtheorem{remark}[theorem]{Remark}

\newcommand{\R}[1][3]{\mathbb{R}^{#1}}
\newcommand{\s}[1][1]{\mathbf{S}^{#1}}

\newcommand{\cross}{\times}
\newcommand{\bndry}{\partial}

\newcommand{\pii}{\pi_1}

\newcommand{\tld}{\tilde}
\newcommand{\br}{\overline}
\newcommand{\tc}{\cal T}

\newcommand{\homeo}{\cong}

\newcommand{\hbdy}{handlebody }

\newcommand{\pfcbc}{$p$-fold cyclic branched cover }

\begin{document}
\title{3-manifolds built from injective handlebodies.}

\author{James Coffey }
\address{University of Melbourne}
\email{coffey@ms.unimelb.edu.au}

\author{Hyam Rubinstein}
\address{University of Melbourne}
\email{H.Rubinstein@ms.unimelb.edu.au}

\begin{abstract}
    This paper looks at a class of closed orientable 3-manifolds constructed from a gluing of three handlebodies, such that the inclusion of each handlebody is $\pii$-injective.  This construction is the generalisation to handlebodies of the condition for gluing three solid tori to produce non-Haken Seifert fibered 3-manifolds with infinite fundamental group. It is shown that there is an efficient algorithm to decide if a gluing of handlebodies meets the disk-condition.   Also an outline for the construction of the characteristic variety (JSJ decomposition) in such manifolds is given.  Some non-Haken and atoroidal examples are given.
\end{abstract}


\maketitle

\section{Introduction}

This paper is concerned with the class of 3-manifolds that meet the disk-condition.  These are closed orientable 3-manifolds constructed from the gluing of three handlebodies, such that the induced map on the fundamental group of each of the handlebodies is injective. Thus all manifolds that meet the disk-condition have infinite fundamental group. The disk condition is an extension to handlebodies of conditions for the gluing of three solid tori to produce non-Haken Seifert fibered manifolds with infinite fundamental group.  These manifolds appear to have many nice properties. In this paper some tools for understanding manifolds that meet the disk-condition are investigated. A number of constructions are given for this class of manifolds including some that are non-Haken and some that are atoroidal.  The characteristic variety of manifolds that meet the disk condition is also investigated.  It is shown that the handlebody structure in fact carries all the information for the characteristic variety.

In section \ref{section:definitions and preliminaries} standard definitions that are used throughout this paper are given.  Also the `disk-condition' is defined and discussed.  In particular it is shown how this condition is a generalisation of the construction of  non-Haken Seifert fibered manifolds with infinite fundamental group.  We also discuss how, on a `gut instinct' level, the class of manifolds that meet the disk-condition will contain many other non-Haken examples.

Section \ref{section:Conditions and examples} is broken up into three subsections.  The first develops some basic tools and also shows that all 3-manifolds that meet the disk-condition have infinite fundamental group and are irreducible.  In the second subsection a sufficient condition is given for gluings of handlebodies to meet the disk-condition. This condition is easily checked and useful for constructing examples.  We then give a necessary and sufficient condition and an algorithm that can be checked in bounded time. The final part is a couple of examples of constructions of manifolds that meet the disk-condition, using Dehn fillings along knots in $\s[3]$ and $n$-fold cyclic branched covers of knots in $\s[3]$.  Some non-Haken examples are produced.

Section \ref{section: char var} is concerned with the characteristic variety in manifolds that meet the disk-condition. The main theorem proved in section \ref{section: char var} is:

\begin{theorem}\label{thrm: torus}
Let $M$ be a closed orientable 3-manifold that meets the disk-condition and $T$ be a torus. If $f:T\to M$ is a $\pi_1$-injective map, then there is $\Sigma \subseteq M$ a Seifert fibered sub-manifold with essential boundary and a map $g:T\to M$ homotopic to $f$ such that  $g(T)\subset \Sigma$.
\end{theorem}

If the characteristic variety has non-empty boundary then the boundary components are essential embedded tori. Therefore a direct corollary of the above theorem is:

\begin{corollary}
If $M$ is a 3-manifold that meets the disk-condition and there is a $\pii$-injective map of the torus into $M$ then either there is a $\pii$-injective embedding of a torus in $M$ or M is a Seifert fibered manifold.
\end{corollary}

These are not new results. However, the aim is to examine how the characteristic variety behaves in manifolds that meet the disk-condition. The proof of the torus theorem is constructive and  gives an algorithm for finding the characteristic variety of manifolds that meet the disk-condition. When the characteristic variety is constructed, the components come in two distinct `flavours'.  The intersection of all three handlebodies in the manifold is a set of injective simple closed curves, called the triple curves.  The first flavour is components which are disjoint from the triple curves.  These components  look very much like the objects that W. Jaco and P. Shalen used to prove the torus theorem for Haken manifolds, see \cite{ja1}.  In each handlebody the components of the characteristic variety are either essential Seifert fibered submanifolds or $I$-bundles.  This is not surprising for if we remove an open neighbourhood of the triple curves we get a manifold with boundary, which is therefore Haken.  Also what is left of the boundaries of the handlebodies is a  set of disjoint spanning surfaces. Therefore the fact these carry all the information for the characteristic variety disjoint from the triple curves is not surprising.

The second flavour of characteristic variety is what we will refer to as the disk components.  In the component handlebodies they look like the regular neighbourhood of intersecting meridian disks.  For this flavour of characteristic variety to occur the manifold must meet a minimal disk-condition, as described in section \ref{section:definitions and preliminaries}. The two flavours of characteristic variety are not necessarily disjoint.  If they do intersect their fiberings can always be made to agree.  In fact, when they intersect, the disk components look like thickened compressing annuli of the components disjoint from the triple curves.

The authors would like to thank Ian Agol for a very helpful comment on this project. We would also like to thank the referee for his/her extremely diligent work which has greatly improved this paper. This research was partially supported by the Australian Research Council.

\section{Definitions and preliminaries}\label{section:definitions and preliminaries}

Throughout this paper we will assume that, unless stated otherwise, we are working in the PL category of manifolds and maps.  Even though we will not explicitly use this structure we will use ideas that are a consequence, such as regular neighbourhoods and transversality as defined by C. Rourke and B. Sanderson in \cite{rou&sa1}. The standard definitions in this field, as given by J. Hempel in \cite{Hem1} or W. Jaco in \cite{ja1}, are used.

A manifold $M$ is \textbf{closed} if $\bndry M = \emptyset$ and \textbf{irreducible} if every embedded $S^2$ bounds a ball.  We will assume, unless otherwise stated that all 3-manifolds are orientable. The reason for this is that all closed non-orientable $\mathbb P^2$-irreducible 3-manifolds are Haken. (A manifold is $\mathbb P^2$-irreducible  if it is irreducible and does not contain any embedded $2$-sided projective planes).

If $M$ is a 3-manifold and $S$ is some surface, which is not a sphere, disk or projective plane,  the map $f:S\to M$ is called $\mathbf{\pi_1}$\textbf{-injective}  if the induced map $f_*:\pii(S)\to \pii(M)$ is injective. If the image of $S$ is not boundary parallel then the map is called an \textbf{essential} map.   Also $f:S \to M$ is a \textbf{proper map} if $f^{-1}(\bndry S) = \bndry S$.  If $F:S\cross I \to M$ is a homotopy/isotopy such that $F(S, 0)$ is a proper map, then it is assumed, unless otherwise stated  that $F(S,t)$ is a proper map for all $t\in I$.  To reduce  notation, an isotopy/homotopy of a surface $S\subset M$ is used without defining the map.  Here we are assuming that there is a map $f:S \to M$ and we are referring to  an isotopy/homotopy of $f$. Defining the map is often unnecessary and would only add to excessive book keeping.

If $H$ is a handlebody and $D$ is a properly embedded disk in $H$ such that $\bndry D$ is essential in $\bndry H$ then $D$ is a  \textbf{meridian disk} of $H$.  If $D$ is a proper singular disk in $H$ such that $\bndry D$ is essential in $\bndry H$, then it is called a \textbf{singular meridian disk}.

In this paper normal curve theory, as defined by S. Matveev in \cite{Mat1}, is used to list finite classes of curves in surfaces.  This definition uses a triangulation of the surface to define normal curves. The surfaces may have polygonal faces, however a barycentric subdivision will produce the required triangulation.

\subsection{The disk-condition}\label{subsection: disk condition}

Before we look at what is meant by the `disk-condition' in closed 3-manifolds, we want to define some objects we need and what is meant by the disk-condition in handlebodies.

\begin{defn}
    For $H$ a handlebody, $\tc$ a set of curves in $\bndry H$ and $D$ a meridian disk, let $|D|$ be the number of intersection between $D$ and $\tc$.
\end{defn}

\begin{defn}
    If $H$ is a handlebody and $\tc$ is a set of essential disjoint simple closed curves in $\bndry H$ then $\tc$ meets the $\mathbf n$ \textbf{disk-condition} in $H$ if for every meridian disk $D$, $|D|\geq n$.
\end{defn}

This seems like a difficult condition to meet, for if $H$ has genus two or higher  there is an infinite number of meridian disks to check.  We later give some sufficient conditions that are easily checked and an algorithm that determines if the disk condition is satisfied.

Next we are going to give a description of the construction of 3-manifolds that meet the `disk-condition'. Please note that even though this description is technically correct it is not enlightening, so later we discuss different ways of looking at these manifolds that are much more useful.

Let $H_1$, $H_2$ and $H_3$ be three handlebodies.  Let $S_{i,j}$, for $i\not= j$ be a sub-surface of $\bndry H_i$ such that:
\begin{enumerate}
    \item
        $\bndry S_{i,j} \not= \emptyset$.
    \item
        The induced map of $\pii(S_{i,j})$ into $\pii(H_i)$ is injective.
    \item
        For $j\not= k$, $S_{i,j}\cup S_{i,k} = \bndry H_i$,
    \item
        $\tc_i  = S_{i,j}\cap S_{i,k} = \bndry S_{i,j} = \bndry S_{i,k}$ is a set of disjoint essential simple closed curves that meet the $n_i$ disk-condition in $H_i$,
    \item
         $S_{i,j}\subset \bndry H_i$ is homeomorphic to $S_{j,i} \subset \bndry H_j$.
\end{enumerate}

\begin{figure}[h]
   \begin{center}
       \includegraphics[width=6cm]{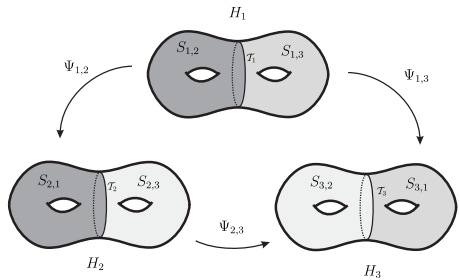}
    \end{center}
   \caption{Homeomorphisms between boundaries of handlebodies.} \label{fig:handlebodies}
\end{figure}

Note that $S_{i,j}$ need not be connected. Now that we have the boundary of each handlebody cut up into $\pii$-injective faces  we want to glue them together by homeomorphisms, $\Psi_{i,j}: S_{i,j}\to S_{j,i}$, that agree along the $\tc_i$'s (figure \ref{fig:handlebodies}).  The result is a closed 3-manifold $M$, for which the image of each handlebody is embedded.

\begin{defn}
    If $M$ is a manifold constructed from three handlebodies as above such that $\tc_i$ meets the $n_i$ disk-condition in $H_i$ and

    \begin{equation}\label{equation: disk-condition}
        \sum_{i=1,2,3} \frac{1}{n_i} \leq \frac{1}{2}
    \end{equation}

    then $M$ meets the $\mathbf{(n_1,n_2,n_3)}$ \textbf{disk-condition}. If we are not talking about a specific $(n_1,n_2,n_3)$, the manifold is said to meet just the \textbf{disk-condition}.
\end{defn}

As we said before, the above definition is not enlightening. Thus, from this point on, we will view 3-manifolds that meet the disk-condition in the following way.  Assume that $M$ is a manifold that meets the disk-condition and $H_1$, $H_2$ and $H_3$ are the images of the handlebodies of the previous definition in $M$.  Then $M = \bigcup_{i=1,2,3} H_i$ and each $H_i$ is embedded in $M$. Then $X = \bigcup_{i=1,2,3} \bndry H_i$ is a 2-complex that cuts $M$ up into handlebodies. As $X$ is constructed by gluing surfaces along their boundaries, it does not meet the meet the usual definition of a 2-complex.  However the surface can be cut up into cells. Also $\tc = \bigcap_{i=1,2,3} H_i$ is a set of essential disjoint simple closed curves in $M$ that meets the $n_i$ disk-condition in $H_i$ where $\sum_{i=1,2,3} 1/n_i\leq 1/2$.

It may seem a bit odd that we are using the same name for the construction of 3-manifolds and the condition on curves in the boundary of handlebodies.  However, the curve condition is the restriction of the condition on compact closed 3-manifolds to each of its  component handlebodies. When we have an equality  in equation \ref{equation: disk-condition}, the result is the three `minimal' cases for the disk-condition. They are; $(6,6,6)$, $(4,8,8)$ or $(4,6,12)$.  These three are of special interest for if a manifold meets the disk-condition, then it meets at least one of these three. Therefore these are the important cases to consider.  It is also worth noting that unlike Heegaard splittings, we can use three handlebodies of different genera.

\begin{figure}[h]
   \begin{center}
       \includegraphics[width=3cm]{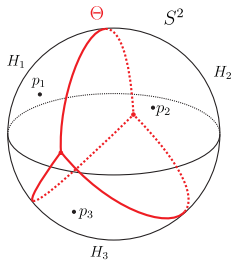}
    \end{center}
   \caption{Base space of non-Haken Seifert fibered space with infinite $\pii$.} \label{fig:Base space of Seifert fibered}
\end{figure}

Another way of viewing a 3-manifold $M$ that meets the disk-condition, is that $X = \bigcup \bndry H_i$ is a 2-complex such that the 1-skeleton $\tc$ consists of essential curves in $X$.  Therefore we can get a manifold $M$ that meets the disk condition by gluing handlebodies to $X$ such that each meridian disk of the handlebodies intersects $\tc$ enough times.  In fact, the disk condition is an extension of the construction of non-Haken Seifert fibered 3-manifolds with infinite fundamental group.  In this case, we know that if a Seifert fibered space is non-Haken with infinite fundamental group, then its base space is a 2-sphere and it has three exceptional fibers of multiplicity $p_i$, where $\sum 1/p_i \leq 1$ (*), as in figure \ref{fig:Base space of Seifert fibered}. This is the construction given by P. Scott in \cite{sc}.  Thus if the inequality (*) is made an equality, the exceptional fibers  have indices $(3,3,3)$, $(2,4,4)$ or $(2,3,6)$.  Another way of viewing this construction is if $\Theta$ is the graph in figure \ref{fig:Base space of Seifert fibered}, then $\Theta \cross \s[1]$ is a 2-complex of three annuli glued together along two triple curves $\tc$.  Then glue in three solid tori $H_i$'s such that the meridian disks wind around $p_i$ times. As we have two triple curves in $\tc$, each meridian intersects $\tc$ $2p_i$ times. Thus, as  $\sum 1/(2p_i)\leq 1/2$, all non-Haken Seifert fibered manifolds with infinite $\pii$ are in the class of manifolds that meet the disk-condition.

Yet another way of viewing 3-manifolds that meet the disk-condition is if we take two handlebodies and glue them together so that we get a 3-manifold with a single incompressible boundary. Then glue a handlebody to this boundary component so the surface is only incompressible in one direction. A very short hierarchy in a closed Haken manifold, as defined by I. Aitchison and H. Rubinstein in \cite{Ai&Ru1}, can be thought of as taking a set of handlebodies, gluing each handlebody to itself so that each of the resulting manifolds has incompressible boundary.  Then glue these incompressible boundaries together to produce the closed manifold. Therefore the incompressible boundaries become the incompressible surface in the Haken manifold. So it would seem that the disk-condition is a weaker condition than being Haken. In fact  we already know this class of manifolds  contains all the non-Haken Seifert fibered manifolds with infinite $\pii$, but it also contains examples of other non-Haken manifolds.

The disk-condition can be easily extended to gluings of four or more handlebodies such that all the statements in this paper follow.  Construct a closed manifold $M$ by gluing together $r \geq 3$ handlebodies, $\{H_1,...,H_r\}$, such that  for $i$, $j$, $k$ and $l$ different; $H_i$ is embedded, $H_i\cap H_j \subset \bndry H_i \cap \bndry H_j$ is a subsurface, $H_i \cap H_j \cap H_k$ is a possibly empty set of pairwise disjoint curves and $H_i \cap H_j \cap H_k \cap H_l = \emptyset$.  Then $X = \bigcup_{1\leq i < j \leq r} H_i \cap H_j$ is a 2-complex which cuts $M$ up into the $H_i$'s and $\tc = \bigcup_{1\leq i < j < k \leq r}  H_i \cap H_j \cap H_k$ is made up of pairwise disjoint simple closed curves. Suppose $\alpha$ is a component of $\tc$. Let $H_{\alpha_1}$, $H_{\alpha_2}$  and $H_{\alpha_3}$ be the three handlebodies around $\alpha$ and $\tc$ meet the $n_{\alpha_i}$ disk-condition in $H_{\alpha_i}$. Then $M$ meets the `generalised disk-condition' if, for each $\alpha \in \tc$, $\sum_{i=1,2,3} 1/n_{\alpha_i} \leq 1/2$. For the purposes of this paper we will not consider such manifolds for $r\geq 4$ as they are all Haken. If $r \geq 4$ then we can  choose $H_i$ and $H_j$ such that $H_i \cap H_j \ne \emptyset$ and there is a component $M'$ of $\overline{M-(H_i \cup H_j)}$ that contains at least two of the handlebodies. Let $S$ be the boundary surface between $H_i\cup H_j$ and $M'$. Then the proof of lemma \ref{lemma: removing pullback graphs} can be altered to show that no closed curve in $S$ bounds a disk and thus $S$  is an embedded incompressible surface.

\section{Conditions and examples}\label{section:Conditions and examples}

For our purposes we need to state Dehn's lemma and the loop theorem in a specific way:

\begin{lemma}  \label{lemma: loop theorem}
    Let $H$ be a handlebody and $\tc$ a collection of essential curves in $\bndry H$. If there is a singular meridian disk $D$ such that $D$ intersects $\tc$ $n$ times then there exists an embedded meridian disk of $H$ that intersects $\tc$ at most $n$ times.
\end{lemma}

Let $H$ be a  handlebody and $\tc$ be a set of disjoint essential simple closed curves in $\bndry H$  that meets the $n$ disk-condition. A direct result of this lemma is that if we have a singular  closed curve $\alpha$ in $\bndry H$ that intersects $\tc$ less than $n$ times and contracts in $H$, then by lemma  \ref{lemma: loop theorem} we know that $\alpha$ is not essential in $\bndry H$.

\begin{lemma}\label{lemma: removing pullback graphs}
    Let $M$ be a manifold that meets the disk condition and $D$ be a disk. If $f:D\to M$ is a map such that $f(\bndry D)\subset int(H_i)$, for some $i$, then $f$ can be homotoped, keeping the boundary fixed, so that $f(D)\subset int(H_i)$.
\end{lemma}

\begin{figure}[h]
   \begin{center}
       \includegraphics[width=6cm]{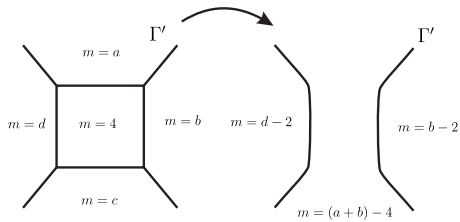}
    \end{center}
   \caption{Removing a $(4,n)$-gon from $\Gamma'$ by homotopy.}
   \label{fig:homotop f to remove gon.}
\end{figure}

\begin{proof}
We will assume that $f(D)$ is transverse to $X$ and that $f(\bndry D)\subset int(H_i)$, where $int(H_i)$ is the interior of the handlebody $H_i$. Thus $\Gamma = f^{-1}(X)$ is a set of trivalent graphs and simple closed curves in $D$. Note that $\bndry D \cap \Gamma = \emptyset$.  Let the $\Gamma_j$'s be the components of $\Gamma$. An innermost component of $\Gamma$, is a component $\Gamma_j$ such that there is a disk $D^*\subset D$ where $\Gamma_j \subset int(D^*)$ and $D^*\cap \Gamma = \Gamma_j$.  Note that if $\Gamma$ is non-empty then it must have at least one innermost component.  Let $\Gamma_j$ be innermost and $D'\subset D^*$ the disk such that $\Gamma_j\subset D'$ and $\bndry D'\subset \Gamma_j$.

If $\Gamma_j$ is a simple loop then $\Gamma_j  = \bndry D'$ and $f(D')\subset H_k$, for $k=1,2$ or $3$. By the disk condition we know that $f(\bndry D')$ must be non-essential  in $\bndry H_k$ as it doesn't intersect $\tc$ and thus $f(D')$ is homotopic into $\bndry H_k$.  We can thus homotop $f$ so that $f(D')\subset \bndry H_k$ and then push it through to remove the component altogether.

If $\Gamma_j$ is a graph then as it is innermost, the faces of $D'$ must all be disks.  Thus each face $F$ of $D'$ is an $(m,n)$-gon, where $F$ has $m$ vertices in its boundary and is mapped by $f$ to a handlebody $H_k$ such that $\tc$ meets the $n$ disk-condition in $H_k$.  We can put a PL metric on $D'$ by assuming that all the edges are geodesic arcs of unit length, that the internal angle at each vertex of an $(m,n)$-gon $F$ is $\pi(1-2/n)$ and all the curvature of $F$ is at a cone point in $int(F)$.   For example if $H$ meets the 6 disk-condition the angle in each corner of an $(m,6)$-gon will be $2\pi/3$.  Note that as each vertex of $\Gamma_i$ in the interior of $D'$ is adjacent to three faces, each is mapped to a different handlebody. Assuming that $M$ meets the $(6,6,6)$, $(4,6,12)$ or $(4,8,8)$ disk-conditions, then the total angle around each vertex in the interior is $2\pi$.  If $F$ is an $(m,n)$-gon, then $\chi(F) = 1$ and the exterior angle sum is $m(2\pi/n)$. If $\mathbf K(F)$ is the curvature of the cone point in $int(F)$, then by the Gauss-Bonnet Theorem we know that:

$$\mathbf K(F) = 2\pi - m(2\pi/n) = 2\pi(1-m/n) $$

Thus for $F$, an $(m,n)$-gon, if $m<n$ then $\mathbf K(F)>0$ or if $m \geq n$ then $\mathbf K(F) \leq 0$. Let $\mathbf F$ be the set of faces of $D'$ and  $\mathbf v$ be the vertices in $\bndry D'$.  For $v\in \mathbf v$ then there are two faces $F_1, F_2 \in \mathbf F$ adjacent to $v$.  Let $F_i$ be an $(m_i,n_i)$-gon. Let the jump angle at $v$ be $\theta_v = \pi - \sum_{i=1,2}\pi(1-2/n_i)$.  By the disk condition $n_i = 4,6,8$ or $12$ and it is not possible to have $n_1 = n_2 = 4$. Thus $\theta_v \leq -\pi/6$. Then once again by Gauss-Bonnet we know that:

$$\sum_{\mathbf F}\mathbf K(F) = 2\pi - \sum_{\mathbf v}\theta_v > 2\pi .$$

This implies that $D'$ must always have some $(m,n)$-gons such that $m<n$.  For example if the manifold meets the $(6,6,6)$ disk-condition then $D'$ would have some $(2,6)$-gons and/or some $(4,6)$-gons.  If $F$ is an $(m,n)$-gon of $D'$ such that  $m< n$ and $f(F)\subset H_k$, then by the disk condition and lemma \ref{lemma: loop theorem}, we know that $f(\bndry F)$ is not essential in $\bndry H_k$.  Thus we can homotop $f$ so that $f(F)$ lies in $\bndry H_k$. We can then homotop $f$ so $f(F)$ is pushed off $\bndry H_k$.  This decreases the total number of faces of $D'$, as shown in figure \ref{fig:homotop f to remove gon.}.   Once again we know that $D'$ has a face of positive curvature that can be removed.  Thus in a finite number of steps $\Gamma_j$ will become a simple closed loop and we can then  homotop $f$ to remove the component $\Gamma_j$ entirely.

As $\Gamma$ always  contains an innermost component, we can continue this process until all of $\Gamma$ has been removed and thus $f(D) \subset int(H_i)$.
\end{proof}

This lemma leads us to a couple of important corollaries about 3-manifolds that meet the disk condition.

\begin{corollary}
    Let $M$ be a 3-manifold that meets the disk condition. Then, for any $1\leq i \leq 3$, the induced map of $\pii(H_i)$ into $\pii(M)$ is injective.
\end{corollary}

\begin{remark}
Note that $\pii(H_i)$ is the free group on $g$ generators, where $g>0$ is the genus of $H_i$. This corollary implies that if a 3-manifold meets the disk-condition, then its fundamental group is infinite.
\end{remark}

\begin{proof}
Let  $D$ be a disk and $\gamma$ be a simple closed curve in $H_i$ that represents a non-trivial element of $\pii(H_i)$.  If the element is trivial in $\pii(M)$, then there would be a map $f:D\to M$ such that $f(\bndry D)=\gamma$.  By lemma \ref{lemma: removing pullback graphs} we can homotop $f$ so that $f(D)\subset int(H_i)$, giving us a contradiction.
\end{proof}

\begin{corollary}
    If $M$ is a 3-manifold that meets the disk-condition, then it is irreducible.
\end{corollary}

\begin{proof}
Let $S$ be a 2-sphere and $f:S\to M$ be an embedding.  Note that $f$ is an embedding and all the moves in the proof of lemma \ref{lemma: removing pullback graphs} can be performed as isotopies. Thus we can isotope $f$ so that $f(S)\cap X = \emptyset$, that is, for some $i$, $f(S)\subset H_i$.  Then, as handlebodies are irreducible, $f(S)$ must bound a 3-ball.
\end{proof}

\subsection{Test for the $n$ disk-condition in handlebodies.}
\label{section: Test for the disk-condition.}

It is not necessary to check every meridian disk of a handlebody to find out if a set of curves in its boundary meets the $n$-disk condition. Let $H$ be a handlebody and $\tc$ the set of curves in $\bndry H$. Let $\mathcal D$ be a set made up of a single representative from each isotopy class of meridian disk of $H$.

The first test is that $\tc$ must separate $\bndry H$ into faces that can be 2-coloured.  Therefore all meridian disks must intersect $\tc$ an even number of times.  From this point on we will assume that $\tc$ is separating in $\bndry H$.

Put a Riemannian metric on $\bndry H$. We will assume that the loops in $\tc$ are length minimizing geodesics. Note that if $\tc$ contains parallel curves, the neighbourhood  of the corresponding length minimizing geodesic can be `flattened', so we can have parallel length minimizing geodesics. We will also assume the boundaries of the disks in $\mathcal D$ are length minimizing geodesics. Both of these can be done simultaneously. From M. Freedman, J. Hass and P. Scott \cite{fhs1} we know that this implies that the number of intersections between any disk in $\mathcal D$ and $\tc$ is minimal, as is the intersection between the boundaries of any two disks in $\mathcal D$. For any disk $D\in \mathcal D$ let $|D|$ be the number of intersections of $\bndry D$ with $\tc$ and for any set of meridian disks $\mathbf D = \{D_i\} \subset \mathcal D$, let $|\mathbf D|= \sum_i |D_i|$.  We can assume that all these curves are transverse to each other. From this point on, unless otherwise stated, when looking at meridian disks we will assume that the number of intersections between their boundaries is minimal.

\begin{lemma}\label{lemma: m-disk intersection.}
    Any two disks of $\mathcal D$ can be isotoped, leaving their boundaries fixed, so that any curves of intersection are properly embedded arcs.
\end{lemma}

\begin{proof}
This proof uses the standard innermost arguments and the fact the handlebodies are irreducible to remove all the components of intersection between two disks that are simple closed curves.
\end{proof}


\begin{defn}
    Let $H$ be a genus $g$ handlebody. We shall call $\mathbf D \subset \mathcal D$ a \textbf{system of meridian disks} if all the disks are disjoint, non-parallel and cut $H$ up into a set of 3-balls. If $\mathbf D$ cuts $\bndry H$ up into $2g-2$ pairs of pants (thrice punctured 2-spheres) then it is a \textbf{basis} for $H$.
\end{defn}

If $H$ has genus $g$, then a minimal system of meridian disks for $H$ consists of $g$ disjoint  non-parallel meridian disks, and the disks cut $H$ up into a single ball.

\begin{defn}
    Let $P$ be a punctured sphere and $\gamma$ be a properly embedded arc in $P$. If  both ends of $\gamma$ are in the one boundary component of $\bndry P$ and the arc is not isotopic into $\bndry P$ then it is called a \textbf{wave}.
\end{defn}

Let $H$ be a handlebody, $\tc$ be a set of essential disjoint simple closed curves in $\bndry H$, $\mathbf D$ be a system of meridian disks for $H$ and $\{ P_1,...,P_l\}$ be the resulting set of punctured spheres produced when we cut $\bndry H$ along $\mathbf D$. Also let $\tc_i = P_i \cap \tc$.  Thus $\tc_i$ is a set of properly embedded disjoint arcs in $P_i$.

\begin{defn}\label{defn: waveless system}
    If each $\tc_i$ contains no waves then $\mathbf D$ is said to be a \textbf{waveless} system of meridian disks for $H$.
\end{defn}

\begin{defn}\label{defn: n-waveless}
    Let $\mathbf D$ be a waveless system of disks. If every wave in each $P_i$ intersects $\tc_i$ at least $n/2$ times, then $\mathbf D$ is called an $\mathbf n$\textbf{-waveless} system of meridian disks.
\end{defn}

If $\mathbf D$ is an $n$-waveless basis then each  $\tc_i$ has at least $n/2$ parallel arcs running between each pair of boundaries in $P_i$.

\begin{lemma}\label{lemma: sufficient disk condition}
    Let $H$ be a handlebody, $\tc\subset \bndry H$ be a separating set of essential simple closed curves and $\mathbf D$ a basis for $H$. If $\mathbf D$ is an $n$-waveless basis, then $\tc$ meets the $n$ disk-condition in $H$.
\end{lemma}

\begin{figure}[h]
   \begin{center}
       \includegraphics[width=3cm]{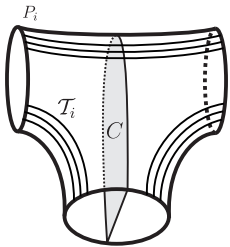}
    \end{center}
   \caption{Bigon in a pair of pants.}
   \label{fig:End-disk in a pair of pants.}
\end{figure}

\begin{figure}[h]
   \begin{center}
       \includegraphics[width=30mm]{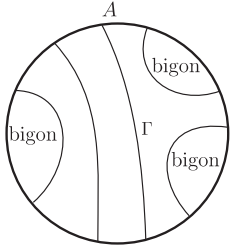}
    \end{center}
   \caption{Meridian disk cut up by arcs of intersection.}
   \label{fig:Meridian disk cut up by arcs of intersection.}
\end{figure}

\begin{proof}
 From the definition of the $n$-waveless condition we know that $\tc$ intersects each disk in $\mathbf D$ at least $n$ times. If $C\in \mathcal D$ is a meridian disk not in $\mathbf D$, then $C \cap \mathbf D \not= \emptyset$.  By lemma \ref{lemma: m-disk intersection.} we can isotop $C$ so that $C \cap \mathbf D$ is a set of  disjoint properly embedded arcs. Therefore if we cut $C$ along $C \cap \mathbf D$ the faces produced must all be disks and contain at least two bigons, as shown in figure \ref{fig:Meridian disk cut up by arcs of intersection.}. Therefore the set $\{ P_i\cap \bndry C\}$ must contain at least two waves, coming from bigons.  As $\mathbf D$ meets the $n$-waveless condition, any wave must intersect $\tc$ at least $n/2$ times, see figure \ref{fig:End-disk in a pair of pants.}. Therefore $\bndry C$ must intersect $\tc$ at least $n$ times.
\end{proof}

If $\tc$ intersects each disk in $\mathbf D$ exactly $n$ times then it must be an  $n$-waveless basis.  The reason is that the only pattern of arcs in a pair of pants, where there are the same number $n$ of endpoints on each boundary curve, consists of $n/2$ arcs joining each pair of boundary loops. This gives us the following corollary.

\begin{corollary}
    Let $H$ be a handlebody, $\tc\subset \bndry H$ be a separating set of simple closed curves and $\mathbf D$ a basis for $H$. If $\tc$ intersects each disk in $\mathbf D$ exactly $n$ times then $\tc$ meets the $n$ disk-condition in $H$.
\end{corollary}

This test for the $n$ disk-condition is a significant restriction. However, it is an easy enough condition to satisfy when constructing examples.

\begin{figure}[h]
   \begin{center}
       \includegraphics[width=8cm]{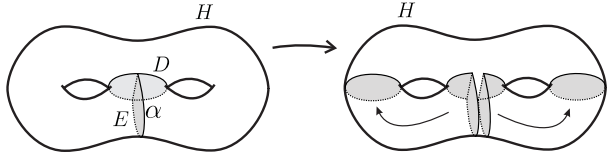}
    \end{center}
   \caption{Boundary compressing a meridian disk.}
   \label{fig:Boundary compressing a meridian disk.}
\end{figure}

Next we describe a specific type of surgery  of  meridian disks. Let $D$ be a meridian disk of $H$ and let $E$ be an embedded disk in $H$ such that $\bndry E \subset D \cup \bndry H$, $\bndry E \cap \bndry D$ is two points, $a_1$ and $a_2$ in $\bndry H$, $\alpha = E \cap \bndry H$ is an arc in  $\bndry H$ which is not homotopic through $\bndry H$ into $\bndry D$ and $D\cap E$ is an arc properly embedded in $D$, as shown in figure \ref{fig:Boundary compressing a meridian disk.}. If we then surger $D$ along $E$ we produce two disks. As $\alpha$ is an arc which is not homotopic through $\bndry H$ into $\bndry D$, both resulting disks are meridian disks isotopic to disks in $\mathcal D$.  We shall call this surgery a \textbf{boundary compression} of a meridian disk.

\begin{figure}[h]
   \begin{center}
       \includegraphics[width=8cm]{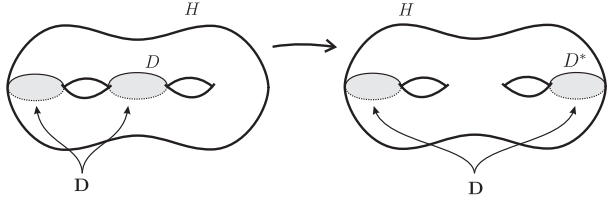}
    \end{center}
   \caption{Disk-swap move.}
   \label{fig:Disk-swap move.}
\end{figure}

Let $\mathbf D$ be a system of disks for the handlebody $H$.  Let $D^*\in \mathcal D$ be a meridian disk disjoint from $\mathbf D$ such that $(\mathbf D \setminus D)\cup D^*$ is a system of meridian disks for some $D\in \mathbf D$.  Then if we remove $D$ from $\mathbf D$ and replace it with $D^*$ this is called a \textbf{disk-swap move} on $\mathbf D$ as shown in figure \ref{fig:Disk-swap move.}.

\begin{lemma} \label{lemma: bndry comp gives m-disk for disk swap move.}
    For a system of meridian disks $\mathbf D =\{D_1, \ldots ,D_n \}$, if we perform a boundary compression on any $D_i$ along a disk disjoint from\\ $\{D_1,... ,D_{i-1},D_{i+1},... ,D_n \}$, then one of the resulting disks can be used for a disk-swap move on $\mathbf D$ removing $D_i$.
\end{lemma}

\begin{remark}
Note that an essential wave in $\overline{\bndry H - \mathbf D}$ defines a  disk swap move on $\mathbf D$.
\end{remark}

\begin{figure}[h]
   \begin{center}
       \includegraphics[width=5.5cm]{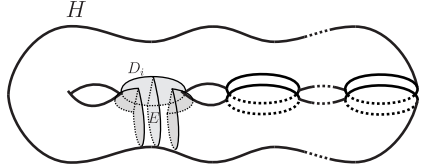}
    \end{center}
   \caption{Boundary compressing a disk from a system of meridian disks.}
   \label{fig:Boundary compressing a disk from a system of meridian disks.}
\end{figure}

\begin{proof}
Let $\mathbf D^*$ be the set of all meridian disks disjoint from $\mathbf D$. Then if a disk $D\in \mathbf D$ is boundary compressed along a disk disjoint from $\mathbf D - D$, one of the resulting disks will be isotopic to a disk in $\mathbf D \cup \mathbf D^*$.   Let $E$ be the disk we are going to boundary compress $D_i$ along.  If we cut $H$ along $\{D_1, \ldots ,D_{i-1},D_{i+1},\ldots ,D_n \}$ the result is a solid torus $T$. Then $D_i$ is a meridian disk of $T$.  Thus a boundary compression on $D_i$ along $E$ will produce two disks, one of which is a meridian disk of $T$ and the other is boundary parallel, as shown in figure \ref{fig:Boundary compressing a disk from a system of meridian disks.}.
\end{proof}

Let $\mathbf D \subset \mathcal D$ be a minimal system of meridian disks for the handlebody $H$. That is, $\mathbf D$ cuts $H$ up into a single ball. Let $\mathbf D^* \subset \mathcal D$ be the set of disks disjoint from $\mathbf D$.

\begin{lemma}\label{lemma: necessary disk condition}
$\tc$ meets the $n$ disk-condition if and only if there is a minimal system of meridian disks $\mathbf D$ such that $|D| \geq n$ for all disks $D\in \mathbf D \cup \mathbf D^*$ and there are no disk swap moves between $\mathbf D$ and $\mathbf D ^*$ that reduce $|\mathbf D|$.
\end{lemma}

\begin{proof}
In the `only if' direction, $\tc$ satisfying the $n$ disk-condition in $H$ implies that $|D|\geq n$ for any $\mathbf D \cup \mathbf D^*$. Given some initial $\mathbf D \cup \mathbf D^*$ we can construct a sequence of disk swaps that reduce $|\mathbf D|$. If $\tc$ meets the $n$ disk condition then such a sequence must terminate, thus giving the required basis.

For the proof in the  `if' direction the first thing to note is that; if there are no disk swap moves to reduce $|\mathbf D|$ then every essential wave in $\overline{\bndry H - \mathbf D}$ must intersect $\tc$ at least $n/2$ times.  Let $D\in \cal D$ be a meridian disk such that $D \not \in \mathbf D \cup \mathbf D^*$.  Then $\Gamma = D\cap \mathbf D \not = \emptyset$.  We are assuming that the intersection between the boundaries of disks is minimal.  Thus by lemma \ref{lemma: m-disk intersection.} we can assume that $\Gamma$ is a set of pairwise disjoint properly embedded arcs in $D$, as shown in figure \ref{fig:Meridian disk cut up by arcs of intersection.}.  $\Gamma$ is minimal and thus all the faces of $D$, when $D$ is cut along $\Gamma$, are disks. Also there must be at least two bigons, $D_1$ and $D_2$. $D_i \cap \overline{\bndry H - \mathbf D}$ are essential waves in $\overline{\bndry H - \mathbf D}$ and thus intersect $\tc $ at least $n/2$ times.
\end{proof}

Next we want to use lemma \ref{lemma: necessary disk condition} to produce an algorithm to determine whether a boundary pattern satisfies the $n$ disk-condition.

\begin{lemma}\label{lemma: waveless minimal system of m-disks}
    Assume we are given a handlebody $H$ and $\tc$ a set of essential curves in $\bndry H$.  There is an algorithm to find, in finite time, a waveless minimal system of meridian disks.
\end{lemma}

\begin{figure}[h]
   \begin{center}
       \includegraphics[width=7cm]{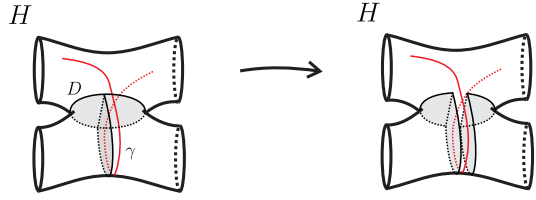}
    \end{center}
   \caption{Boundary compression to remove a wave.}
   \label{fig:boundary_comp_to_remove_wave}
\end{figure}

\begin{proof}
Suppose we start with an arbitrary minimal system of meridian disks $\mathbf D$ for $H$. If $\tc$ has a wave when $H$ is cut along $\mathbf D$, then there is a sub-arc $\gamma\subset \tc$ with both ends in some disk $D\in \mathbf D$ and $int(\gamma)\cap \mathbf D = \emptyset$.  Then $D$ has a boundary compression disk $E$ such that the arc $E \cap \bndry H =\gamma$.  Let $D_1$ and $D_2$ be the disks produced by compressing $D$ along $E$. Then $|D_i| \leq |D| - 2$, as shown in figure \ref{fig:boundary_comp_to_remove_wave}. Thus when a disk swap move is done swapping $D$ for one of the $D_i$'s, $|\mathbf D|$ will decrease by at least two. Note also that the number of waves does not go up. If there is another wave we can always do another boundary disk compression and a disk swap move to reduce $|\mathbf D|$, thus this process must terminate in a finite number of moves.
\end{proof}

\begin{figure}[h]
   \begin{center}
       \includegraphics[width=3cm]{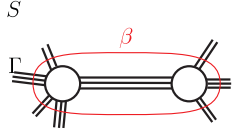}
    \end{center}
   \caption{Boundary of meridian disk to add to $\mathbf D$ .}
   \label{fig:cutting off pants}
\end{figure}

Given that it is possible to find a waveless minimal system of meridian disks $\mathbf D$, to show that we can  find a waveless basis, we proceed as follows.  Suppose we have already found a waveless system of disks and want to add new waveless disks, until we get a basis. We can use our initial set of boundary curves of disks to cut $\bndry H$ to obtain a punctured sphere $S = \overline{\bndry H- \mathbf D}$. Suppose that there is at least one pair of boundary curves of $S$ such that all the arcs of $\Gamma$ running between them are parallel.  Then there is a simple closed curve $\beta$ which it is essential in $S$, is not boundary parallel and each curve in $\Gamma$ intersects $\beta$ at most once, as shown in figure \ref{fig:cutting off pants}. Then we can add a disk with boundary $\beta$ to enlarge our system of waveless disks.

To simplify this problem, collapse each boundary component of $S$ to a vertex and identify parallel copies of edges of $\Gamma$. This produces a graph $\Gamma'$ embedded in  a 2-sphere $S'$ such that $\Gamma'$ is connected, no two edges are parallel and no edge has both ends at the one vertex.  This means that if we cut $S$ along $\Gamma'$ all the resulting faces will be disks and will have degree at least 3.

\begin{figure}[h]
   \begin{center}
       \includegraphics[width=5.5cm]{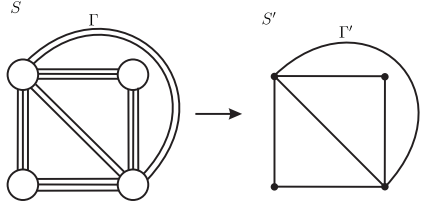}
    \end{center}
   \caption{$\Gamma$ and $\Gamma'$.}
   \label{fig: triple curve to graph.}
\end{figure}

\begin{defn}\label{defn: 2-cycle}
Let a \textbf{2-cycle} in a graph be a simple closed loop that is the union of two edges.
\end{defn}

The problem has now become to show that we can always find two vertices of $\Gamma'$ that are joined by exactly one edge. That is, to find a vertex not contained in a 2-cycle.  Let $c$ be a 2-cycle in $\Gamma'$, thus $c$ cuts $S'$ into two disks and as $\Gamma'$ does not contain any parallel edges, the interior of both disks must contain at least one vertex of $\Gamma'$. We now want to show that there is  a vertex of $\Gamma'$ that is not part of a 2-cycle.  Let $c$ and $c'$ be two 2-cycles in $\Gamma'$.  If $c\cap c'$ is empty, a single vertex or edge, then the interior of one of the disks produced when we cut $S'$ along $c$ must be disjoint from  $c'$.  If $c\cap c'$ is two vertices, then we can construct a third 2-cycle $c''$ such that when we cut $S'$ along $c''$, the interior of one of the disks produced is disjoint from both $c$ and $c'$. If $C$ is the set of all 2-cycles in $\Gamma'$, then there must be a 2-cycle $c\in C$ such that when $S'$ is cut along $c$ we get a disk $D$ such that there are no 2-cycles intersecting $int(D)$.  As there are no parallel edges in $\Gamma'$, $\Gamma' \cap int(D) \ne \emptyset$. Therefore $\Gamma'$ has to have a vertex in $int(D)$ that is not in a 2-cycle.  This gives us the following lemma.

\begin{lemma}\label{lemma: waveless minimal basis of m-disks}
    Assume we are given a handlebody $H$ and $\tc$ a set of essential curves in $\bndry H$.  There is an algorithm to find, in finite time, a waveless basis.
\end{lemma}

Note that this means that once the minimal waveless system of meridian disks has been found, most of the work has been done and that to produce a waveless basis, suitable meridian disks are added to the system.  This lemma is not expressly used in the rest of this paper, but waveless bases are used in section \ref{section: char var} in a condition for atoroidal manifolds. Thus it is nice to know that given a 3-manifold that meets the disk-condition, we can always find a waveless basis for each of its handlebodies.

\begin{lemma}\label{lemma: algorithm to test necessary condition}
    Let $H$ be a handlebody and $\tc$ a set of essential curves in $\bndry H$.  Then there is an algorithm to determine, in finite time, if $\tc$ satisfies the $n$ disk-condition.
\end{lemma}

\begin{proof}
Once again let $\mathbf D$ be a minimal system of disks and $n(\mathbf D)$ be a regular neighbourhood of $\mathbf D$. Let $S = \overline{\bndry H- n(\mathbf D)}$ and $\Gamma = \tc \cap \overline{\bndry H- n(\mathbf D)}$. Then $S$ is a $2g$-punctured sphere, where $g$ is the genus of $H$. Also $\Gamma$ is a set of arcs properly embedded in $S$. By lemma \ref{lemma: waveless minimal system of m-disks} we can assume that $\Gamma$ does not contain any waves.  Therefore $\Gamma$ cuts $S$ up into embedded polygons of degree at least four.  As above let $\mathbf D^* \subset \mathcal D$ be the set of meridian disks disjoint from $D$. For any $D^*\in \mathbf D^*$, $ D^* \cap S = \alpha$ is a simple closed curve in $int(S)$.  Let $|\alpha|$ be the number of times that $\alpha$ intersects $\Gamma$. Note that $|\alpha|=|D^*|$.  We have therefore reduced the question of looking for meridian disks disjoint from $\mathbf D$ to looking at essential simple closed curves in $S$. For $D \in \mathbf D$ then $n(D) \cap S$ is two boundary curves, $\bndry D_1$ and $\bndry D_2$, of $S$.  Then if $\gamma$ is  an essential curve in $S$ that separates $\bndry D_1$ from $\bndry D_2$, the disk bounded by $\gamma$ can be used for a disk swap move on $D$.  Let $N=max \{ |D|:D\in \mathbf D\}$ and $L$  be the set of essential simple closed curves in $S$ of length at most $N$.  Thus as $L$ is a finite set of curves and as each face of $S$ is a polygon, we can list all the elements of $L$ using normal curve theory. Therefore to test whether $\mathbf D$ satisfies lemma \ref{lemma: necessary disk condition}  we need to check that; all disks in $\mathbf D$ intersect $\tc$ at least $n$ times, all the  curves in $L$ have length at least $n$, and $|\gamma|\geq |D|$ for $\gamma \in L$ and $D \in \mathbf D$  such that $\gamma$ separates the two disks $D \cap S$ in $S$. If a disk swap move is found, then we perform the move and then test the new system.  As $|D|$ decreases by at least two with each move, the algorithm will terminate in finite time, either when a suitable system is found, meaning $\tc$ satisfies the $n$ disk-condition or when a meridian disk is found that  intersects $\tc$ less than $n$ times.
\end{proof}

Note that this algorithm can be continued until a system is found which has a `locally minimal' intersection. If $n=min \{|D|:D\in \mathbf D\}$, then $n$ is the supremum disk-condition satisfied by $\tc$.  For if there is a meridian disk that intersects $\tc$ less than $n$ times that is not in $\mathbf D$, then the algorithm would not have terminated.  An equivalent statement is that $\mathbf D$ is an $n$-waveless  system of disks.  Clearly if there is an essential wave in $\overline{\bndry H - \mathbf D}$ that intersects $\tc$ less than $n/2$ times then there is a disk swap move to reduce $|\mathbf D|$. In the other direction, if $\mathbf D$ is an $n$-waveless system and there is a meridian disk $D \in \cal D$ such that $|D|<n$, then clearly $D\cap \mathbf D \ne \emptyset$.  Thus $D$ can be surgered to give a boundary compressing disk for some disk in $\mathbf D$ and thus a wave in $\overline{\bndry H - \mathbf D}$, that intersects $\tc$ at less than $n/2$ points. Therefore there is an alternative algorithm to test the disk condition, giving the corollary:

\begin{corollary}\label{corollary: finding n-waveless system}
    If $H$ is a handlebody and $\tc \subset \bndry H$ is a set of essential curves that meet the $n$ disk-condition, then there is an algorithm to find an $n$-waveless minimal system of meridian disks.
\end{corollary}


\subsection{Examples}

To construct manifolds that meet the disk-condition a technique such as Dehn surgery or branched covers must build a manifold which contains a 2-complex that cuts it up into three injective handlebodies.

\subsubsection{Dehn filling examples.}\label{subsection: dehn filling example}

\begin{figure}[h]
   \begin{center}
       \includegraphics[width=4cm]{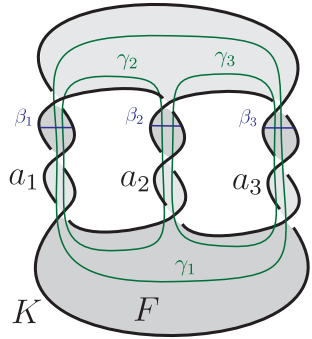}
    \end{center}
   \caption{$(3,3,3)$-pretzel knot.}
   \label{fig:(333)pretzel_knot.}
\end{figure}

\begin{figure}[h]
   \begin{center}
       \includegraphics[width=4cm]{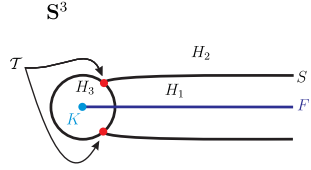}
    \end{center}
   \caption{Handlebodies in Dehn filling construction.}
   \label{fig: Handlebodies in Dehn filling construction.}
\end{figure}

The first class of examples of manifolds that meet the disk-condition are constructed by performing Dehn surgery along suitable knots in $\s[3]$.  Let $K\subset \s[3]$ be the $(3,3,3)$-pretzel knot  and $F$ the free spanning surface shown in figure \ref{fig:(333)pretzel_knot.}.  For $A\subset \s[3]$ let $n(A)$ be a regular neighbourhood of $A$.  Let $H_3=n(K)$ and $H_1 = \overline{n(F) - H_3}$, as in figure  \ref{fig: Handlebodies in Dehn filling construction.}.  Then $H_1$ is a genus 2 \hbdy and $\tc = \bndry (H_1 \cap H_3)$ is two copies of $K$.  $H_1$  is homeomorphic to an $I$-bundle over $F$ and $\tc$ to the boundary curves of the vertical boundary of the $I$-bundle structure.  The three arcs, $\beta_i$'s, in figure \ref{fig:(333)pretzel_knot.} lift to a basis of meridian disks for the I-bundle.  Each wave in the pairs of pants produced when $\bndry H_1$ is cut along the basis intersects $\tc$ at least twice. Therefore the basis is 4-waveless and by lemma \ref{lemma: sufficient disk condition} $\tc$ meets the 4 disk-condition in $H_1$. $H_2=\overline{\s[3] - (H_1\cup H_3)}$ is a genus 2 \hbdy and the three curves, $\gamma_i$'s, in figure \ref{fig:(333)pretzel_knot.} bound meridian disks of a basis $\mathbf D$ for $H_2$. As $\tc$ is two copies of $K$ each wave in the two pairs of pants, produced by cutting $\bndry H_2$ along the $\gamma_i$'s, intersect $\tc$ six times.  Thus $\mathbf D$ is a 12-waveless basis for $H_2$ and by lemma \ref{lemma: sufficient disk condition} $\tc$ meets the 12 disk-condition in $H_2$.  Therefore if a Dehn surgery along $K$ is performed such that the meridian disk of the solid torus glued back in intersects $\tc$ at least 6 times, a manifold that meets the $(4,6,12)$ disk-condition is produced.   U. Oertel showed in \cite{Oe1} that all but finitely many Dehn surgeries on such pretzel knots produce non-Haken 3-manifolds.

This construction can  be generalised to any knot  $K\subset \s[3]$, that has a free spanning surface $F$, such that $K$ meets the 6 disk-condition in $\overline{\s[3] - F}$.  Then any Dehn surgery of type $(p,q)$ with $|p|\geq 6$ will produce a manifold meeting the disk-condition.


\subsubsection{Branched cover examples.}

The next method for constructing manifolds which meet the disk-condition is taking cyclic branched covers over knots in $\s[3]$.  We look at two conditions on knots that are sufficient for the resulting manifolds to meet the disk condition.

Let $B_i$'s for $i=1,2$ or 3 be 3-balls and $\gamma_i = \{\gamma_i^1,...,\gamma_i^k\}$, for $k\geq 2$, be a set of properly unknotted pairwise disjoint embedded arcs in $B_i$.  Unknotted means that there is a set of pairwise  disjoint embedded disks, $D_i = \{D_i^1,...,D_i^k\}$,  such that $\gamma_i^j \subset  \bndry D_i^j$ and $\overline{\bndry D_i^j - \gamma_i^j} = D_i^j \cap \bndry B$.  Therefore,  if we take the $p$-fold cyclic branched cover of $B_i$, with $\gamma_i$ as the branch set, then the result will be a genus $(p-1)(k-1)$ handlebody $H_i$.  Let $r_i:H_i \to B_i$ be the branched covering map and $\alpha_i \subset \bndry B_i$ be a simple closed loop disjoint from $\gamma_i$ such that $\tc_i = r^{-1}(\alpha_i)$ meets the $n_i$ disk-condition in $H_i$.  Note that $\alpha_i$ can be thought of as cutting $\bndry D_i$ up into two hemispheres.

\begin{figure}[h]
   \begin{center}
       \includegraphics[width=6cm]{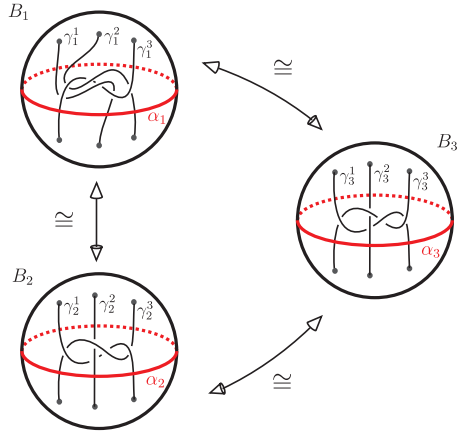}
    \end{center}
   \caption{Bubble construction.}
   \label{fig:bubble construction.}
\end{figure}

Now glue the three balls by homeomorphisms between their hemispheres, as shown in figure \ref{fig:bubble construction.}, so that the resulting manifold is $\s[3]$ and the endpoints of $\gamma_i$'s match up.  Thus $K = \bigcup \gamma_i$ is a link and $C =\bigcup \bndry B_i$ is a 2-complex of three disks glued along a triple curve $\alpha$, which is the image of the $a_i's$.  Let $M$ be the $p$-fold cyclic branched cover of $\s[3]$ with $K$ as the branch set.  Let $r:M\to \s[3]$ be the branched covering map.  Then $X=r^{-1}(C)$ is a 2-complex that cuts $M$ up into handlebodies and $\tc = r^{-1}(\alpha)$ is a set of triple curves that meets the $n_i$ disk condition in $H_i$. Thus if $\sum \frac{1}{n_i} \leq \frac{1}{2}$, $M$ meets the disk-condition.

\begin{figure}[h]
   \begin{center}
       \includegraphics[width=4cm]{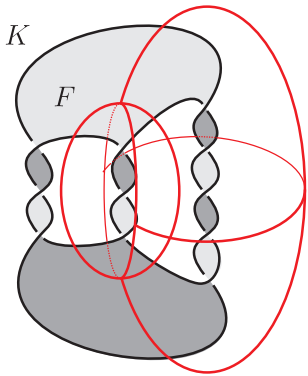}
    \end{center}
   \caption{$(3,3,5)$-pretzel knot.}
   \label{fig:(335)pretzel_knot.}
\end{figure}

If  $k=2$ or $3$ and the intersection of $\alpha_i$ with $D_i$ is minimal under isotopy in $\bndry B_i - \gamma_i$,   then a sufficient condition for the lift of $\gamma_i$ to the \pfcbc  of $B_i$ to meet the $n$ disk-condition is that any essential wave in $\overline{\bndry B_i - D_i}$ intersects $\gamma \cap \overline{\bndry B_i - D_i}$ at least $n/2$ times.  Note that this is a slight variation of lemma \ref{lemma: sufficient disk condition} and the proof is essentially the same.  Given the 2-complex shown in figure  \ref{fig:(335)pretzel_knot.}, it can be seen that any \pfcbc over an $(a_1,a_2,a_3)$-pretzel knot in $\s[3]$ such that $|a_i|\geq2$, will produce a manifold that meets the disk-condition.

Let $M$ be a manifold that meets the disk-condition and can be constructed from the gluing of three genus 2 injective handlebodies. Then a simple Euler characteristic argument shows that all the faces of the 2-complex $X$ must either be once punctured tori or twice punctured disks. If all the faces are once punctured tori then the set of triple curves, $\tc$, is a single curve. Thus a $\pi$ involution of  $\tc$ can be canonically extended, up to isotopy, to each of the faces of $X$. This is done via a waveless basis for each handlebody. Thus the involution can be extended to the whole of $M$.  This means that any such manifold has  a $\mathbb Z_2$ symmetry and is the 2-fold cyclic branched  cover of $\s[3]$ over some link. Also if all the faces of $X$ are pairs of pants then there is an involution of $M$ where the fixed point set is a graph of with nine edges and three vertices, where each vertex is order three. The quotient space is once again $\s[3]$.

The second construction is done by taking the 3-fold branched covers of the knots that meet essentially the same condition as in the Dehn filling construction and then the lift of the Seifert surface gives the 2-complex.  Let $K$ be a knot in $\s[3]$ and $F$ be a free Seifert surface of $K$. That means that $\overline{\s[3] - F}$ is a handlebody.  For these examples we use the construction of 3-fold cyclic branched covers over knots in $\s[3]$ given by D. Rolfsen in \cite{ro1}.  Let $n(K)$ be a regular neighbourhood of $K$, $\alpha \subset \bndry n(K) $ be the meridian curve of $n(K)$ and $N = \overline{ \s[3] - n(K)}$.  Let $\tld N$ be the 3-fold cyclic cover of $N$ and $p: \tld N \to N$ the covering projection.  That is, let $G \subset \pii(N)$ be the kernel of the homomorphism mapping $\pii(N)$ onto $\mathbb Z_3$, where the meridian of $n(K)$ is sent to the generator.  Then $\tld N$ is the cover corresponding to $G$.  This means that  $\tld N$ has a single torus boundary and $\tld \alpha = p^{-1} (\alpha)$ is a single curve that covers $\alpha$ three times.  Therefore $\tld F= p^{-1}(F)$ is a set of three properly embedded spanning surfaces in $\tld N$.  As $F$ is free,  $\overline{\tld N - \tld F}$ is three handlebodies. Let $M$ be the 3-fold cyclic branched cover of $\s[3]$ with $K$ as the branch set.  Then $M$ can be constructed by gluing a solid torus $T$ to $\bndry \tld N$ so that its meridian matches $\tld \alpha$.   Next extend each surface in $\tld F$ along an annulus to the spine $\tc$ of $T$ to produce a 2-complex $X$.  Thus $X$ is a 2-complex that cuts $M$ into three handlebodies. Thus for $M$ to meet the disk-condition it is sufficient for $K$ to meet the 6 disk-condition in $\overline{\s[3] - F}$. An obvious example of such a knot is the $(3,3,3)$ pretzel knot in figure  \ref{fig:(333)pretzel_knot.}.

The 3-fold cyclic branched cover of the $(3,3,5)$ pretzel knot $K$ is an example of a manifold with two distinct splitting 2-complexes that meet the disk-condition.  Let $M$ be the  3-fold cyclic branched cover of $\s[3]$ with $K$ as the branch set. Let $X$ be the 2-complex produced by lifting the Seifert surface $F$ to $M$ and let $X'$ be the 2-complex produced by lifting the `bubble' 2-complex shown in figure  \ref{fig:(335)pretzel_knot.}.   $X$ and $X'$ are distinct 2-complexes meeting the disk-condition. That is there is no homeomorphism of $M$ that sends $X$ to $X'$, for if there was, $M$ would have a $\mathbb Z_3 \cross \mathbb Z_3$ symmetry and thus $K$ would have a $\mathbb Z_3$ symmetry, which is clearly not the case.  Note that if each twisted band in $K$ has the same number of crossings, for example the $(3,3,3)$ pretzel knot, then the 3-fold cyclic branched cover does have a $\mathbb Z_3 \cross \mathbb Z_3$ symmetry.

\section{Characteristic Variety} \label{section: char var}

In this section we  prove the torus theorem and construct the characteristic variety in 3-manifolds that meet the disk-condition. The first step is to look at how, in the component handlebodies,  properly embedded essential annuli disjoint from the triple curves intersect and how meridian disks that intersect $\tc$ exactly $n_i$ times intersect.  This allows us to build a picture of the characteristic variety in each of the handlebodies, which we then use to construct the characteristic variety of the manifold.

\subsection[Embedded annuli and meridian disks.]{Handlebodies, embedded annuli and meridian disks. }

Throughout  this section let $H$ be a handlebody and $\tc$ be a set of disjoint essential simple closed  curves in $\bndry H$ that meet the $n$ disk-condition in $H$. Also when an annulus $A$ is said to be properly embedded in a handlebody $H$, it is assumed that it is disjoint from $\tc$. We will also assume that all intersections between surfaces are transverse. Before we look at the components of the characteristic variety in each handlebody, we need to look at some properties of embedded essential annuli that are disjoint from the triple curves.

\subsection{Essential annuli.}

In this section we are looking at some properties of intersections between embedded essential proper annuli.

\begin{defn}
An intersection curve between two annuli is said to be \textbf{vertical} if it is a properly embedded arc which is not boundary parallel in both annuli. The intersection curve is \textbf{horizontal} if it is a simple closed essential loop in both annuli.
\end{defn}

\begin{figure}[h]
\begin{center}
   \subfigure[Horizontal]{
     \includegraphics[width=3cm]{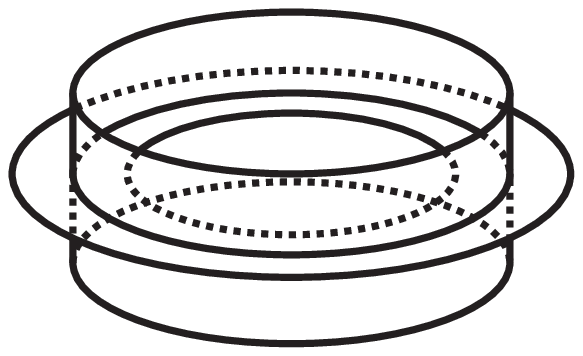}}
   \hspace{1cm}
   \subfigure[Vertical]{
     \includegraphics[width=4.5cm]{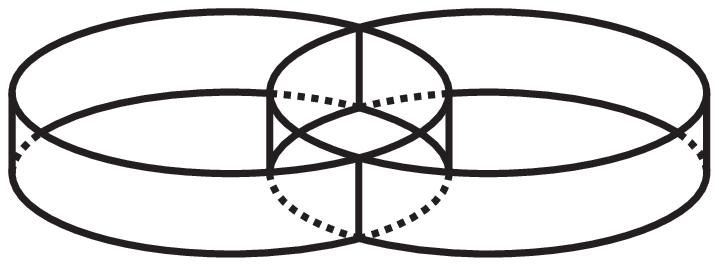}}
\end{center}
\caption{Intersecting embedded annuli.}
\label{fig:Intersecting embedded annuli.}
\end{figure}

If there is a proper isotopy in $H-\tc$ of two annuli which removes their intersection, then the annuli  will be referred to as having \textbf{trivial intersection} and if the intersection cannot be removed, the annuli have \textbf{non-trivial intersection}.  This means that if two embedded annuli have non-trivial intersection they cannot be isotopically parallel. The disk-condition restricts how properly embedded annuli can intersect.


\begin{lemma}\label{lemma: intersecting annuli}
Let $A_1$ and $A_2$ be two essential properly embedded annuli in $H$.  Then there is a proper isotopy of them in $H-\tc$ so that all their intersections are either vertical or horizontal.
\end{lemma}

\begin{remark}
This means that non-trivial intersections between embedded annuli must either be all horizontal or all vertical.
\end{remark}

\begin{proof}
This proof is done by constructing the isotopy using the usual innermost arguments and the following observations.  Let $A_1$ and $A_2$ be essential properly embedded  annuli in $H - \tc$ and let $\Gamma = A_1\cap A_2$. If the intersections between the $A_i$'s are all horizontal, then all the components are simple closed loops and essential in both $A_i$'s and if all  the intersections are vertical then all the components of $\Gamma$ are arcs that run between the boundary curves of the $A_i$'s.   First note that as the $A_i$'s are embedded they cannot have both horizontal and vertical intersections. As $H$ is irreducible there is an isotopy of $A_1$ to remove components of $\Gamma$ that are simple closed loops and inessential in both $A_i$'s.  Also by irreducibility of $H$ and the disk-condition, there is an isotopy of $A_1$ to remove components  of $\Gamma$ which are properly embedded arcs and boundary parallel in both $A_i$'s.  Let $\gamma$ be a component of $\Gamma$ such that it is a simple closed loop and essential in $A_1$ and not essential in $A_2$. Then the disk in $A_2$ bounded by $\gamma$ gives a compression of $A_1$.  As the resulting disks are disjoint from $\tc$, they must be parallel into $\bndry H$ and thus $A_1$ is not essential in $H$. Now let $\gamma$ be a component of $\Gamma$ such that it is a properly embedded arc that has both ends in the same boundary curve of $A_1$ and runs between the boundary curves of $A_2$.  Then the disk produced by cutting $A_1$ along $\gamma$ is a boundary compression disk for $A_2$ and the disk produced by compressing $A_2$ is disjoint from $\tc$, thus implying that $A_2$ is boundary parallel in $H -\tc$.
\end{proof}


\begin{lemma}\label{lemma: embedded annuli int vert or horz}
Let $H$ be a handlebody and $\tc$ a  set of curves of $\bndry H$ that meet the $n$ disk-condition. Assume a  properly embedded essential annulus intersects two  other properly embedded essential annuli, one vertically and the other horizontally. Then if there is a non-trivial horizontal intersection, the vertical intersections can be removed by an isotopy.
\end{lemma}

\begin{remark}
This indicates there are three types of essential embedded annuli in $H$;  those that have non-trivial horizontal intersections with other annuli, those that have non-trivial vertical intersections with other annuli and those that have no non-trivial intersections with other annuli. Later in this section we will see that these types of annuli correspond to the flavours of characteristic variety in $H-\tc$.
\end{remark}

\begin{figure}[h]
  \begin{center}
      \includegraphics[width=5.5cm]{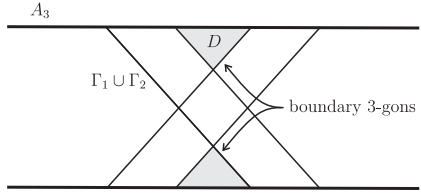}
   \end{center}
  \caption{Component of the pull back graph $\Gamma_1 \cup \Gamma_2$.}
  \label{fig: Component of the pull back graph.}
\end{figure}

\begin{proof}
Let $A_1$, $A_2$ be two properly embedded essential annuli that have non-trivial horizontal intersection. Let $A_3$ be the third embedded essential annulus that intersects $A_1$ vertically. If we  assume that the  vertical intersection between $A_1$ and $A_3$ is non-empty,  $(A_1\cap A_2) \cap A_3 \not = \emptyset$ and thus the intersection between $A_2$ and $A_3$ is non-empty.  By lemma   \ref{lemma: intersecting annuli} we can isotope this intersection to be either vertical or horizontal. If the intersections are horizontal, $\bndry A_3$ is disjoint from $\bndry A_2$ as both $A_2\cap A_1$ and $A_2\cap A_3$ are essential simple closed curves in $A_2$. There is an innermost bigon on $A_2$ bounded by one arc from each of $A_2\cap A_1$ and $A_2\cap A_3$ with common endpoints, see figure \ref{fig: Curves of intersection in A_2.}. This is clear because each arc of $A_1\cap A_3$ has to have at least one corresponding vertex of $(A_2\cap A_1)\cap(A_2\cap A_3)$. It is then straightforward to see that there are vertical arcs of intersection of $A_1\cap A_3$ which contain the two vertices of this bigon, see figure \ref{fig: Curves of intersection in A_3.}. We can then isotope $A_3$ along this bigon to convert these two vertical arcs into two boundary parallel arcs of $A_1\cap A_3$ which can be removed by a further isotopy. In this way, eventually all the vertical arcs of $A_1\cap A_3$ can be removed. Thus we can assume that $A_3$ intersects both $A_1$ and $A_2$ vertically.

\begin{figure}[h]
  \begin{center}
      \includegraphics[width=3cm]{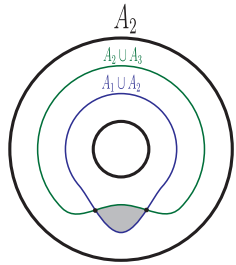}
   \end{center}
  \caption{Curves of intersection in $A_2$.}
  \label{fig: Curves of intersection in A_2.}
\end{figure}

\begin{figure}[h]
  \begin{center}
      \includegraphics[width=7cm]{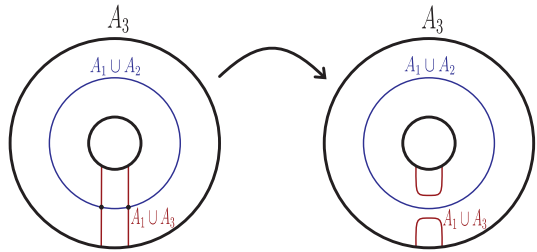}
   \end{center}
  \caption{Curves of intersection in $A_3$.}
  \label{fig: Curves of intersection in A_3.}
\end{figure}

Let  $\Gamma_i = A_3 \cap A_i$, for $i\not= 3$. Then $\Gamma_i$ is a set of properly embedded pairwise disjoint spanning arcs in $A_3$, where each arc from $\Gamma_1$ intersects at least one arc from $\Gamma_2$.  The faces produced when $A_3$ is cut up  along $\Gamma_1 \cup \Gamma_2$ are all disks.  As each connected component of $\Gamma_1 \cup \Gamma_2$ contains at least two arcs, each component will have a boundary 3-gon, $D$ , as shown in figure \ref{fig: Component of the pull back graph.}, such that sub-arcs of  $\bndry A_3$, $\Gamma_1$ and $\Gamma_2$ make up its three edges.  Then the disk $D$ gives an isotopy of $A_1$ that converts the corresponding essential closed curve of $A_1\cap A_2$ into a boundary parallel arc. Thus there is a further isotopy to remove the intersection altogether.  This process can be repeated to remove all the intersections of $A_1\cap A_2$, giving a contradiction.
\end{proof}

Therefore if a proper essential annulus in $H - \tc$ has a non-trivial horizontal/vertical intersection with one annulus, then we can arrange that all its non-trivial intersections with all other essential annuli must be horizontal/vertical.

\subsection{Meridian disks}

Next we want to look at some properties of intersecting meridian disks.  In particular if $\tc$ meets the $n$ disk-condition in $H$, then there may be meridian disks that intersect $\tc$ exactly $n$ times.  These disks are important when we are dealing with the disk flavour of characteristic variety.

\begin{defn}
If $F$ is an $n$-gon and $\gamma$ is a properly embedded arc in $F$ such that if $F$ is cut along $\gamma$, the result is two disks that intersect $\tc$ $n/2$ times, then $\gamma$ is said to be a \textbf{bisecting} arc of $F$.
\end{defn}

\begin{lemma}\label{lemma: intersection between n-gons bisecting}
Let $H$ be a handlebody and $\tc$ a set of curves in $\bndry H$ that meets the $n$ disk-condition.  If $D_1$ and $D_2$ are meridian disks that intersect $\tc$ $n$ times, then there is an isotopy of $D_1$ and $D_2$ so that $\Gamma = D_1 \cap D_2$ is a set of  properly embedded disjoint bisecting arcs in both $D_i$'s or the intersection can be removed.
\end{lemma}

\begin{figure}[h]
  \begin{center}
      \includegraphics[width=3cm]{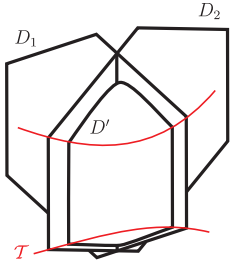}
   \end{center}
  \caption{Two trivially intersecting 6-gons.}
  \label{fig:2_intersecting_6-gons.}
\end{figure}

\begin{proof}
This proof uses the usual innermost arguments and the following observations, to construct an isotopy to remove arcs of $\Gamma$ that are not bisecting in both disks. By lemma \ref{lemma: m-disk intersection.} we can assume that all components of $\Gamma$ are properly embedded arcs.  Let $\gamma$ be an innermost arc of $\Gamma$ which is not bisecting in $D_1$.  Let $D$ be the disk produced by cutting $D_1$ along $\gamma$ such that $D$ intersects $\tc$ less than $n$ times. Then one of the disks produced by surgering $D_2$ along $D$ must intersect $\tc$ in less than $n$ points, as shown in figure \ref{fig:2_intersecting_6-gons.}, and thus is boundary parallel in $H$.  So there is an isotopy of $D_1$  to remove $\gamma$.
\end{proof}

\begin{lemma}\label{lemma: no triple points in intersecting n-gons}
Let $H$ be a handlebody, $\tc$ be a set of curves in $\bndry H$ that meets the $n$ disk-condition and $D_1$, $D_2$ and  $D_3$ a set of meridian disks that all intersect $\tc$ $n$ times.   Then there is an isotopy of the  $D_i$'s so that $\bigcap D_i = \emptyset$.
\end{lemma}

\begin{proof}
By the previous lemma we can isotope $D_1$ and $D_2$ so that their intersection is a set of parallel bisecting arcs in both disks.  Assume that all the trivial intersections between $D_1$ and $D_2$  have been removed and that  $D_1\cap D_2 \not = \emptyset$.  Let $A$ be the regular neighbourhood of $D_1\cup D_2$ and $B$ be the frontier of $A$ in $H$.  As any annulus of $B$ does not intersect $\tc$, $B$ is a set of meridian disks that intersect $\tc$ exactly $n$ times and essential annuli whose boundary compressing disks intersect $\tc$ at least $n/2$ times.

Let $D$ be a disk and $f:D\to H$ be an embedding such that $f(D) = D_3$.  Then $\Gamma=f^{-1}(B)$ is a set of properly embedded pairwise disjoint curves.  As usual there is a homotopy of $f$ to remove components of $\Gamma$ that are simple closed curves.  If $D_3$ intersects an annulus of $B$ then, from above, either the intersections are bisecting parallel arcs or there is a homotopy of $f$ to remove them. Similarly from lemma \ref{lemma: intersection between n-gons bisecting} if $D_3$ intersects a disk of $B$ then either the intersections are bisecting parallel arcs or there is a homotopy of $f$ to remove them.  Therefore there is a homotopy of $f$ so that $\Gamma$ is a set of parallel bisecting arcs and thus $f^{-1}(D_2\cup D_3)$ is a set of parallel bisecting arcs.
\end{proof}

\subsection[Flavours of characteristic variety.]{Flavours of characteristic variety in the handlebodies.}

\subsubsection{$I$-bundle regions}

Let $H$ be a handlebody and $\tc$ a set of essential simple closed curves in $\bndry H$, that meet the $n$ disk-condition in $H$.  Let $N$ be a maximal, up to isotopy, I-bundle in $H$ disjoint from $\tc$, with its horizontal boundaries embedded in $\bndry H -\tc$ and each component of $N$ has non-trivial fundamental group and the induced map on the fundamental group is injective. Thus $N$ is an I-bundle with a base space which is an embedded surface in $H$. Let $S$ be a component of this embedded surface. If $S$ is orientable then the corresponding component  of $N$ has a product structure and its horizontal surface consists of two copies of $S$ embedded in $\bndry H - \tc$. Alternatively if $S$ is non-orientable then the corresponding component of $N$ has a horizontal boundary which is a double cover of $S$ embedded in $\bndry H - \tc$. In both cases the vertical boundary is a set of essential properly embedded annuli. From this point on they will be called boundary annuli. Also note that none of the base surfaces can be disks. This means that $N$ is a set of embedded handlebodies in $H$ with genus $\geq 1$. $N$ is not unique, for if $H$ contains two embedded annuli that intersect horizontally, in a non-trivial way, then $N$ can contain the regular neighbourhood of one or the other annulus but not both.

\begin{defn}
Let the  $\mathbf I$\textbf{-bundle region}, $N_I$, be the set of all components $N_i$ from $N$ which have base spaces that are not annuli or Mobius bands.
\end{defn}

Later the $I$-bundle region is shown to be unique up to isotopy.


\begin{lemma}\label{lemma:embedded vert int annuli into I-bundle}
If $A$ is a properly embedded essential in $H - \tc$ that has a non-trivial vertical intersection with another properly embedded essential annulus, then it is isotopic into $N_I$.
\end{lemma}

\begin{proof}
Let $A$ be an annulus and the map $f_i:A \to H -\tc$, for $i= 1$ or $2$ be an essential proper embedding such that $f_1(A)=(A_1)$ and $f_2(A)=A_2$ have  non-trivial vertical intersections. Let $B$ be the set of boundary annuli of $N_I$. If $A_1\cap N_I \ne \emptyset$ then by lemma \ref{lemma: intersecting annuli}  and lemma \ref{lemma: embedded annuli int vert or horz} we know that there is an isotopy  of $f_1$ so that the intersection between $A_1$ and the annuli in $B$ is vertical.  Thus the pullback $\Gamma_1= f_1^{-1}(B)$ is a set of properly embedded non-boundary parallel curves in $A$ and as $B$ is separating in $H$, there must be an even number of them.   Thus $\Gamma_1$ cuts $A$ up into quadrilaterals and every alternate one is mapped by $f_1$ into $\overline{(H-N_I)}$.  Let $A'\subset A$ be a quadrilateral such that $f(A')\subset \overline{(H-N_I)}$.  Also let $n(f_1(A'))$ be the regular neighbourhood of $f_1(A')$ in $\overline{(H-N_I)}$ disjoint from $\tc$.  Note that $n(f_1(A'))$ can be fibered as an $I$-bundle over a quadrilateral.  Then there must be an isotopy of $f_1$ to remove the curves $\Gamma_1 \cap A'$ otherwise $n(f_1(A'))\cup N_I$ would be   larger than $N_I$, contradicting maximality.  We can repeat this process until $\Gamma_1 = \emptyset$, thus $A_1\cap B = \emptyset$. This process can be repeated for $A_2$ so that it is disjoint from $B$.  If $A_1\cap A_2$ is disjoint from $N_I$ then $n(A_1\cup A_2)$ can be fibered as an $I$-bundle and added to $N_I$, contradicting maximality.  Thus $A_1 \cup A_2 \subset N_I$.
\end{proof}

Note that the above lemma suggests that if an annulus $A$ meets another horizontally, it may not be possible to isotop $A$ into $N_I$.


Now let $\breve H$ be a regular finite sheeted cover of $H$ and $\breve{\tc}$ be the lift of $\tc$.  Thus $\breve H$ also is a handlebody with $\breve{\tc}$ satisfying the $n$ disk-condition.  Now let $N_I \subset \breve H$ be the $I$-bundle region, as described above.  Also let $G$ be the group of covering translations of $\breve H$ such that $\breve H/G = H$.  Let $N_i$, for $1\leq i\leq n$, be the connected sub-handlebodies of $N_I$ and $S_i$ be the base-surface corresponding to $N_i$.

\begin{lemma} \label{lemma: I-bundle in regular cover}
If  $N_i$ is a component of $N_I$, then for any $g\in G$, $g(N_i)$ is isotopic to a component of $N_I$.
\end{lemma}

From the previous lemma we get the following corollary.

\begin{corollary}\label{corollary:  g(N_I) is isotopic to N_I}
For any $g\in G$, $g(N_I)$ is isotopic to $N_I$.
\end{corollary}

This corollary can be used to show that $N_I$ can be isotoped so that it is preserved by $G$. Put a Riemannian metric  on $H$, lift it to $\breve H$ and then isotope $N_I$ so that the boundary annuli of the $N_I$ are least area. Let $g\in G$ and $A$ be a boundary annulus of $N_I$. By the arguments used by  Freedman, Hass and Scott in \cite{fhs2}, $g(A)$ is either a boundary annulus of $N_I$ or disjoint from all boundary annuli of $N_I$. Let $N_I'$ and $N_I''$ be components of $N_I$ such that $g(N_I')$ is isotopic  to $N_I''$. If $N_I'\not= N_I''$, then replace $N_I''$ by $g(N_I')$.  Now assume that $N_I' =N_I''$.   We need to look at what happens to the boundary annuli under $g$.  Let $A$ and $A'$ be boundary annuli of $N_I'$ such that $g(A)$ is isotopic to $A'$.  If $A\not= A'$ then replace $A'$ by $g(A)$. Now assume that $A=A'$ and $g(A)\not= A$.   As each element of $g$ is a homeomorphism, $g(N_I')\not\subset int(N_I')$.  Then by this observation and maximality  of $N_I$,  $g(N_I')\cap N_I'$ is empty or isotopic to $N_I'$. Another way of saying this is that $\overline{g(N_I') - N_I'}$ and $\overline{N_I' - g(N_I')}$ are sets of thickened annuli.  We can then assume that $g(A)$ is disjoint from $N_I'$.  Let $U_i$, for $i\in \mathbb N$, be the thickened annulus component of $\overline{g^i(N_I') - g^{i-1}(N_I')}$, where $g^0$ is the identity.  As $\breve H$ is a finite sheeted normal cover, there is some $m\in \mathbb N$ such that $g^m$ is the identity.  Therefore $U_1 \cup ...\cup U_m$ is an  annulus bundle over $\s[1]$ properly embedded in $\breve H$, which cannot happen, thus $g(A) = A$. This gives us the following corollary.

\begin{corollary}\label{corollary: N_I preserved under cover transformations}
There is an isotopy of $N_I \subset \breve H$ such that it is preserved by all the covering transformations.
\end{corollary}

Lemma \ref{lemma:embedded vert int annuli into I-bundle}  implies that if $H$ contains two embedded annuli that have non-trivial vertical  intersection then $N_I$ is not empty.  Note this is a sufficient condition not a necessary one.  For example if $N_I$ is an $I$-bundle over a twice punctured disk then any embedded annuli contained in $N_I$ are parallel to boundary annuli  and thus their intersections can be removed.

\begin{proof}\textit{(of lemma \ref{lemma: I-bundle in regular cover})}
Let $\mathbf{A}$ be the set of boundary annuli of $g(N_i)$ and $\mathbf{B}$ be the set of boundary annuli of $N_I$.  If  $g(N_i)$ and $N_I$  have a non-trivial intersection, then by lemma \ref{lemma: intersecting annuli} there is an isotopy of $g$ so that if any annuli in $\mathbf{A}$ and any annuli in $\mathbf{B}$ intersect, then the intersection curves are all either vertical or horizontal.  Now isotope $g$ to remove all trivial intersections between annuli in $\mathbf A$ and $\mathbf B$.

Let $B\in \mathbf{B}$ be an annulus such that it  intersects at least one annulus in $\mathbf{A}$ horizontally.  By lemma \ref{lemma: embedded annuli int vert or horz} it can only intersect the other annuli in $\mathbf{A}$ horizontally.  Thus $B\cap g(N_i)$ is a set of annuli properly embedded in $g(N_i)$. Let $B'\subset B$ be one such annulus.

Isotop $B'$ so that it is transverse to the $I$-bundle structure. As intersections of $B$ with annuli in $\mathbf A$ are minimal, $B'$ either projects 1-to-1 onto the base space or double covers it.  This depends on whether the two boundary curves of $B'$ are in different annuli in $A$ or in the same annulus, respectively. Therefore the base space of $g(N_i)$ and thus $N_i$ is either an annulus or a Mobius band, giving us a contradiction. This means that all horizontal intersections between annuli in $\mathbf A$ and  $\mathbf B$ can be removed.

Therefore all intersections between annuli in $\mathbf A$ and  $\mathbf B$ that are non-trivial are vertical.  But by lemma \ref{lemma: singular annulus} we can isotope all such annuli in $\mathbf A$ into $N_I$. Therefore there is an isotopy of $g$ so that $g(N_i)\cap N_I \neq \emptyset $ and $\mathbf A \cap \mathbf B = \emptyset $.  Thus we know that we can isotope $g$ so that $g(N_i)$ lies inside $N_I$,  otherwise $g(N_i) \cup N_I$ would be a larger $I$-bundle than $N_I$, contradicting maximality.

As $g(N_i)$ is connected we know that it lies in a single components, $N_k$, of $N_I$. If $g(N_i)$ is not isotopic to $N_k$ then $g^{-1}(N_k-g(N_i)) \cup N_I$ is a larger $I$-bundle region, contradicting maximality.
\end{proof}

\begin{lemma}
$N_{I}$ is unique up to ambient isotopy of $H$.
\end{lemma}

We will not do the proof for this lemma as the technique is the same as lemma \ref{lemma:  I-bundle in regular cover}.  The idea being that if we  assume that we have two  $I$-bundle regions $N_I$ and $N_I'$ that are not isotopic then we get a contradiction to their maximality. Another property of $N_I$ we need later is the lemma:

\begin{lemma}\label{lemma: compressing disk of i-bundle has int at least n/2}
Let $H$ be a handlebody, $\tc$ be a set of pairwise disjoint essential simple closed curves in $\bndry H$ that meet the $n$ disk-condition and $N_I$ be the  $I$-bundle region in $H$. Then if $A$ is a boundary annulus of $N_I$ and $D$ is a boundary compression disk for $A$ then $|D|\geq n/2$.
\end{lemma}

\begin{figure}[h]
  \begin{center}
      \includegraphics[width=6cm]{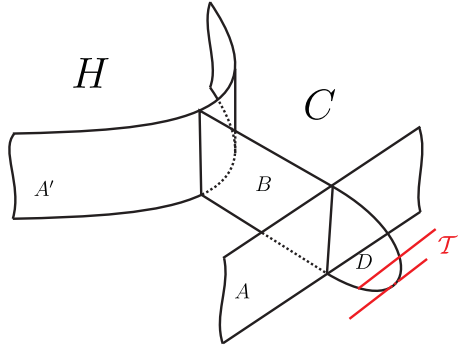}
   \end{center}
  \caption{Extending boundary compression disk through an $I$-bundle component.}
  \label{fig: Extending boundary compression disk.}
\end{figure}

\begin{proof}
Assume that $N_I$ has a boundary annulus $A$ with a boundary compressing disk $D$ such that $|D|<n/2$. Also let $N_i$ be the component of $N_I$ that has $A$ as a boundary annulus.  If we compress $A$ along $D$ to get a disk $E$ then $|E|<n$.  Therefore $A$ must be boundary parallel, meaning there is a proper isotopy of $A$ into $\bndry H$. Note that this does not mean there is a proper isotopy of $A$ into $\bndry H - \tc$. First assume that $N_i$ has more than one boundary annulus.  Let $A'$ be another boundary annulus of $N_i$. As $N_i$ is an $I$-bundle there is a 4-gon $B$ properly embedded in $N_i$, such that $B\cap A = D\cap A$ and $A'\cap B$ is a properly embedded arc in $A'$ that is not boundary parallel, as shown in figure \ref{fig: Extending boundary compression disk.}, for suitable choice of $D$.  Let  $D'=D\cap B$. Then $|D'|\leq n/2$  and if we compress $A'$ along $D'$ we get a disk isotopic to $E$. Therefore $A$ and $A'$ must be parallel and $N_i$ is the regular neighbourhood of a properly embedded annulus and thus can not be contained in $N_I$.  If $N_i$ has a single boundary annulus $A$, then similarly by the $I$-bundle structure, there is a properly embedded 4-gon $B \subset N_I$ such that it is not boundary parallel and $A\cap B$ is two arcs that are not parallel into $\bndry A$.  Then there are two boundary compression disks for $A$ that can be glued to $B$ along $A \cap B$. This produces a meridian disk that intersects  $\tc$ less than $n$ times, contradicting the disk-condition.
\end{proof}

\subsubsection{Tree Regions}

Now let $N=\{N_i\}$ be a maximal set, up to isotopy, of  fibered solid tori embedded in $H- \tc$ such that $N_i\cap N_j = \emptyset$ for $i\neq j$ and $\bndry H  \cap N_i$ is a non-empty set of annuli that are essential in both $\bndry N_i$ and $\bndry (H) - \tc$ and $\overline{\mathrm{int} (H) \cap \bndry N_i}$ is a non-empty set of annuli not isotopic into $\bndry H - \tc$ for each $i$.  $N$ is the maximal tree region of $H-\tc$.  This name will become clearer when we describe it further. Note that by Haken-Kneser finiteness arguments we can see that $N$ has a finite number of components.

\begin{defn}
Let a \textbf{simple }$\mathbf q$\textbf{-tree} be a tree that is the cone on $q$ points.  A vertex of valency one is called an \textbf{end vertex}.
\end{defn}

Let $Q$ be a simple $q$-tree. Embed $Q$ in $\R[2]\subset \R[3]$.  Let $P^Q$ be a $2q$ polygon embedded in $\R[2]$ such that every alternate edge intersects $Q$ at an end vertex.  Colour the edges of $P^Q$ containing an end vertex of $Q$ green and all the others red.  Then let $A_q = P^Q\cross [0,1]$ and $a_t=P^Q\cross\{t\}$, for $t =0$ or $1$. Let $\Phi _p$ be a homeomorphism between $a_0$ and $a_1$ that twists by $\frac{2\pi}{p}$, such that it maps green edges to green edges and red to red.  This means that $p = \frac{q}{n}$ for $n \in \mathbb Z$. Let $A_{(p,q)}$ be $A_q$ with the faces $a_0$ and $a_1$ glued according to $\Phi_p$.  Therefore $A_{(p,q)}$ is a torus fibered by $\s$ with an exceptional fiber of order $(p,q)$. For each $N_i\in N$ there is a unique $(p_i,q_i)$ such that there is a fiber preserving homeomorphism from $A_{(p_i,q_i)}$ to $N_i$ where the fibering agrees with the boundary curves of the boundary annuli.

\begin{figure}[h]
  \begin{center}
      \includegraphics[width=6cm]{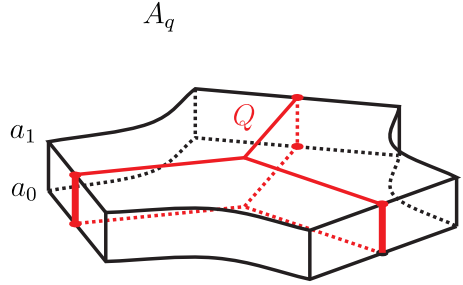}
   \end{center}
  \caption{An example of an $A_q$.}
  \label{fig: An example of a A_q.}
\end{figure}

Let $A_1$ and $A_2$ be two properly embedded essential annuli in $H-\tc$ that intersect horizontally and $n(A_1\cup A_2)$ be a regular neighbourhood disjoint from $\tc$.  Then $\overline{\bndry n(A_1\cup A_2) \cap int H}$ is a set of properly embedded annuli and tori.  Let $T$ be such a torus. The induced map on $\pii(T)$ is not trivial and $\pii(H)$ does not contain any free abelian subgroups of rank 2. Therefore  $T$ is unknotted  and bounds a solid torus whose intersection with $n(A_1\cup A_2)$ is $T$.  Glue solid tori to each torus boundary  of $\overline{\bndry n(A\cup B) \cap int H}$ to produce a submanifold $P$.  Now $\overline{\bndry P \cap int H}$ is a set of properly embedded essential annuli and $P$ is a solid torus.   Note  there is a homeomorphism from $P$ to some $A_{(p,q)}$ that sends the boundary curves of  $P \cap \bndry H$ to fibers of $A_{(p,q)}$.

\begin{defn}
 Let the \textbf{tree region},  $N_T$,  be the union of all components  $N_i\in N$ such that $p_i > 2$.
\end{defn}

As with the $I$-bundle region, we are removing the components of $N$ that are homeomorphic to $A_{(1,2)}$ or $A_{(2,2)}$, that is regular neighbourhoods of properly embedded annuli or Mobius bands, to get $N_T$. This is because if there are two annuli in $H-\tc$ that have a non-trivial vertical intersection then the maximal tree region can contain the regular neighbourhood of only one of the annuli.  Therefore the maximal tree region is not unique. Later it is shown that the tree region is unique up to isotopy.


\begin{lemma}\label{lemma:embedded horz int annuli into tree region}
If $A$ is a properly embedded annulus in $H$ that has at least one non-trivial horizontal  intersection with another properly embedded annulus, then there is an isotopy of $A$  into $N_T$.
\end{lemma}

This proof is similar to lemma \ref{lemma:embedded vert int annuli into I-bundle}.

\begin{proof}
Let $A$ be an annulus and the map $f_i:A \to H$, for $i= 1$ or $2$ be an essential proper embedding such that $f_i(A)=A_i$ is disjoint from $\tc$ for each $i$ and $A_1$ and $A_2$ have  non-trivial horizontal intersections. Let $B$ be the set of boundary annuli of $N_T$. If $A_1\cap N_T \ne \emptyset$ then by lemma \ref{lemma: intersecting annuli}  and lemma \ref{lemma: embedded annuli int vert or horz} we know that there is an isotopy  of $f_1$ so that the intersection curves between $A_1$ and the annuli in $B$ are horizontal. Thus the pullback $\Gamma_1= f_1^{-1}(B)$ is a set of  essential simple closed curves in $A$. Therefore $\Gamma_1$ cuts $A$ up into essential annuli.  Let $A'\subset A$ be one of these annuli such that $f_1(A')\subset \overline{H-N_T}$ and let $n(f_1(A'))$ be the regular neighbourhood of $f(A')$ disjoint from $\tc$. Then $n(f_1(A'))$ can be fibered as an $A_{(1,2)}$ fibered torus.  Thus there must be an isotopy of $f_1$ to remove the curves $ A'\cap \Gamma_1$ ( there may be just one if $\bndry A \cap \bndry A' \ne \emptyset$) otherwise $N_T\cup n(f_1(A'))$ would be larger than $N_T$, contradicting maximality. So by repeating this process, there is an isotopy of $f_1$ such that $ A_1 \cap B = \emptyset$.  This same process produces an isotopy of $f_2$ so that $A_2 \cap B = \emptyset$.  If $A_1\cup A_2$ is disjoint from $N_T$ then, as above, the torus boundaries of $n(A_1\cup A_2)$ can be filled in with solid tori so the resulting manifold $P$ is a solid torus. Then $N_T \cup P$ will be a larger tree region contradicting maximality, thus $A_1 \cup A_2 \subset N_T$.
\end{proof}


Once again let $\breve H$ be a finite sheeted normal cover of $H$, $\breve{\tc}$ be the lift of $\tc$ and $G$ be the group of covering translations of $\breve H$ such that $\breve H / G = H$. Also let $N_T$ be the tree region in $\breve H$. We then get the following lemma.

\begin{lemma}\label{lemma: tree region in regular cover}
Let $N_i$ be a component of $N_T$. For any $g\in G$, $g(N_i)$ is isotopic to an element of $N_T$.
\end{lemma}

From  the previous lemma we get the following corollary.

\begin{corollary}\label{corollary: g(N_T) isotopic to N_T}
For any $g\in G$, $g(N_T)$ is isotopic to $N_T$.
\end{corollary}

From the above corollary and using the same least area arguments as we did with $I$-bundle regions we get the following corollary.

\begin{corollary}
There is an isotopy of $N_T$ in $\breve H$ so that it is preserved by the covering transformations.
\end{corollary}

This means that $N_T$ will project down to a non-trivial tree region in $H$. If $H$ contains two embedded annuli that have a non-trivial horizontal intersection then $H$ has a non-empty reduced tree region.  Note this is a sufficient condition but not a necessary one. The following is similar to the proof of lemma \ref{lemma: I-bundle in regular cover}.

\begin{proof}\emph{(of lemma \ref{lemma: tree region in regular cover})}
Assume that $N_i$ is a component of $N_T$ and for some $g\in G$, $g(N_i)$ is not isotopic to an element of $N_T$.  Let $\mathbf{A}$ be the set of boundary annuli of $g(N_i)$ and $\mathbf{B}$ be the set of boundary annuli of $N_T$. By lemma \ref{lemma: intersecting annuli} we know that there is an isotopy of $g$ so that any annuli from $\mathbf A$ and $\mathbf B$ intersect vertically or horizontally.  Also remove all trivial intersections.

Let $B$ be an annulus in $\mathbf B$ that intersects annuli from $\mathbf A$ vertically.  Then $B\cap g(N_i)$ is a set of properly embedded squares in $g(N_i)$.  Let $B'$ be one such square. As the number of intersections between $B$ and $\mathbf A$ has been minimized $\bndry B'$ is essential in $\bndry g(N_i)$.  Therefore $g(N_i)$ and thus $N_i$ is either the regular neighbourhood of an annulus or Mobius band.  This implies that $p_i = 2$, contradicting that $N_i$ is a component of $N_T$. Then any intersections between annuli from $\mathbf A$ and $\mathbf B$  must be non-trivial and horizontal. By lemma \ref{lemma:embedded horz int annuli into tree region} we can isotope all such annuli from $\mathbf A$ into $N_T$.

We have now isotoped $g$  so that $\mathbf A\cap \mathbf B= \emptyset$.  We can thus isotope $g$ so that $g(N_i)$ lies inside a single component of $N_T$, otherwise $g(N_i)\cup N_T$ would be a larger tree region, contradicting maximality of $N_T$.  Let $g(N_i)$ lie in $N_k\in N_T$. If $g(N_i)$ is not isotopic to $N_k$ then $g^{-1}(N_k-g(N_i)) \cup N_T$ is a larger tree region.
\end{proof}

\begin{lemma}
$N_{T}$ is unique up to ambient isotopy of $H$.
\end{lemma}

We will not do the proof for this lemma as the working is the same as lemma \ref{lemma:  I-bundle in regular cover}.  The idea being that if we  assume that we have two tree regions $N_T$ and $N_T'$ that are not isotopic then we get a contradiction to their maximality.


\subsubsection{Annulus regions}
\label{section: Annulus region}

It is clear from the definitions of $N_I$ and $N_T$ that:

\begin{lemma}
If $H$ is a handlebody and $\tc$ is a set of curves in $\bndry H$ that meet the $n$ disk-condition, then there is an isotopy of $N_I$ and $N_T$ so that $N_I \cap N_T = \emptyset$.
\end{lemma}

Let $A_I$ be the set of $I$-bundles in a maximal $I$-bundle region but not in $N_I$. That is, they have base spaces that are either annuli or Mobius bands.  Let $A_T$ be the set of fibered  tori that are in a maximal tree region but not in $N_T$. This is, they are all the components of the maximal tree region whose associated trees have two end vertices. Then let $N_A$ be the set of components from $A_T$ and $A_I$ that are isotopically equivalent in $H-\tc$. Components of $N_A$ are regular neighbourhoods of properly embedded annuli or Mobius bands and that they can be fibered by intervals or circles.  The components of $A_I - N_A$ ($A_T - N_A$) are the components of the maximal $I$-bundle (maximal tree region) that cause the maximal $I$-bundle (maximal tree region) to be not unique and, in fact, the components of $A_I - N_A$ ($A_T - N_A$) can be isotoped into $N_T$ ($N_I$).

Clearly by the definition, $N_A$ can be isotoped to be disjoint from $N_I$ and $N_T$.  Therefore it is contained in the set of handlebodies $H'= \overline{H -(N_I\cup N_T)}$.  Any annulus that can be made to intersect another non-parallel annulus either vertically or horizontally is isotopic into  $N_I\cup N_T$.  Thus any non-parallel annuli in $H'$ cannot be isotoped to intersect either vertically or horizontally.  Therefore by the maximality of the maximal $I$-bundle region and the maximal tree region we know that $N_A$ is isotopic to the regular neighbourhood of the maximal set of properly embedded annuli in $H'$. Thus we get the following lemma.

\begin{lemma}
 $N_A$ is unique up to ambient isotopy  of $H$ and can be isotoped to be disjoint from $N_I\cup N_T$.
\end{lemma}

\begin{defn}
If $H$ is a handlebody and $\tc$ is a set of triple curves in its boundary that meets the $n$ disk-condition, then for the pair $\{H,\tc\}$ let  the  \textbf{maximal annulus region}  be $N=N_I\cup N_T \cup N_A$ where $N_I$, $N_T$ and $N_A$ are as defined above.
\end{defn}



\subsubsection{Disk Regions}

In this section we want to define the building blocks for the flavour of characteristic variety that intersects the triple curves.   In each handlebody $H_i$ they look like the regular neighbourhood of meridian disks that intersect the  triple curves exactly $n_i$ times,  where $\sum 1/n_i = 1/2$.  Hence we will refer to them as \textbf{disk regions}.  Let $H$ be a handlebody and $\tc$ a set of essential curves in its boundary that meet the $n$ disk-condition in $H$. Let $\mathbf D$ be a set made up of a single representative from each isotopy class of meridian disks that intersect $\tc$ exactly $n$ times. Assume that the disks in $\mathbf D$ have been isotoped so that the intersection between any pair of disks is a set of bisecting arcs and the intersection between any three disks is empty. Let $n(\mathbf D)$ be the regular neighbourhood of $\mathbf D$.  Then $\overline{\bndry n(\mathbf D) \cap int(H)}$ is a set of properly embedded disks that intersect $\tc$ $n$ times and annuli that are disjoint from $\tc$. For any of the boundary components that are either non-meridian disks or non-essential annuli, add the appropriate 3-cell to $n(\mathbf D)$. The resulting sub-manifold $P$ is the  \textbf{disk region}.

\begin{figure}[h]
  \begin{center}
      \includegraphics[width=6cm]{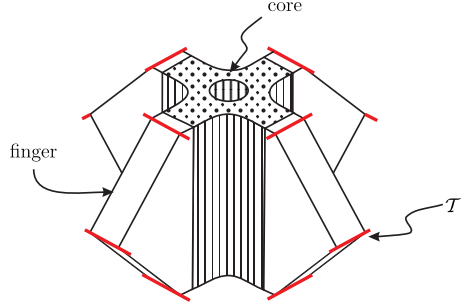}
   \end{center}
  \caption{A component of the disk region.}
  \label{fig:2_disk region}
\end{figure}

By lemma \ref{lemma: no triple points in intersecting n-gons} we can isotope the disks in $\mathbf D$ so that the intersection between any pair of disks is a set of parallel bisecting arcs and the intersection between any three is empty.  Therefore, for any disk $D_i \in \mathbf D$, $\Gamma_i = D_i \cap (\mathbf D\backslash D_i) $ is a set of parallel bisecting arcs.  Let $D_i'$ be the disk produced when $D_i$ is cut along its innermost arcs.  Let $\mathbf D'$ be the set of disks produced when this is done to all disks in $\mathbf D$.  Then $\bigcup D_i'$ is an $I$-bundle over a graph.  This fibering can then be extended to the 'core' of each component of $P$.  The un-fibered parts of each component are the regular neighbourhood of disks that intersect $\tc$ $n/2$ times and which boundary compress  the boundary annuli of the core.  We will call these fingers, see figure \ref{fig:2_disk region}.  Note that each component has at least one finger.  Unlike the $I$-bundle regions defined earlier, the core may have a disk as its base space.  The fibering of each component is unique, up to isotopy, except if the component is the regular neighbourhood of a single meridian disk. In the latter case we do not fiber the core until later.

\begin{lemma}\label{lemma: singular disks into C}
All singular meridian disks which intersect $\tc$ $n$ times can be homotoped into $P$.
\end{lemma}

\begin{proof}
Let $D$ be a disk and $f:D\to H$ be a singular map such that $A=f(D)$ is a singular meridian disk. Let $P$ be the maximal disk region, as defined above and $f^{-1}(\tc)$ be $n$ vertices in $\bndry D$.  Then $B = \overline{\bndry P \cap int(H)}$ is a set of order $n$ meridian disks and  annuli essential in $H-\tc$.  Then $\Gamma = f^{-1}(B)$ is a set of properly embedded arcs and simple closed curves in $D$. As $H$ is irreducible there is a homotopy of $f$ to remove all simple closed curves from $\Gamma$.  Thus $\Gamma$ is a set of properly embedded disjoint simple arcs in $D$.

By maximality of $P$ any boundary compressing disks of a component of $B$, as described in section \ref{section: Test for the disk-condition.}, must intersect $\tc$ more than $n/2$ times.  There must be an innermost disk $D_1\subset D$ such that $f(D_1)$ intersects $\tc$ at most $n/2$ times.  Thus by Dehn's lemma and the loop theorem, see lemma \ref{lemma: loop theorem}, we can remove any arcs from $\Gamma$.  We can repeat this process until  $A$ is disjoint from $B$.  Thus either $A$ is contained in $P$ or disjoint from $P$. If it is disjoint then there must be an isotopy of $f$ so that $A\subset P$. Otherwise using Dehn's lemma and the loop theorem  we get a contradiction to the maximality of $P$.
\end{proof}

\subsection{Handlebodies and singular annuli}
\label{section: Handlebodies and immersed annuli}

In  W. Jaco and P. Shalen's proof of the torus theorem, an essential step is the annulus theorem. In fact, the torus theorem is a consequence of the annulus theorem. Similarly, a lemma that is a slight variation of the annulus theorem is required here. This variation is simplified as it is restricted to handlebodies. Namely, suppose a handlebody $H$ has a set of curves in its boundary, $\tc$, that meets the $n$ disk-condition. Assume also there is a proper essential ``singular'' map $f$ of an annulus into $H-\tc$. Then $f$ is homotopic to a proper essential ``singular'' map of an annulus into  the maximal annulus region.  There are two main steps to prove this lemma.  The first is to show that if there is a proper singular essential map of an annulus into $H-\tc$ then there is a similar embedded one.   Next we show any proper essential embedding of an annulus in  $H - \tc$ is  isotopic to one in its maximal annulus region.


\begin{lemma}\label{lemma: singular annulus}
Let $H$ be a handlebody and $\tc$ a set of simple closed curves in $\bndry H$ that meet the $n$ disk-condition. Let $A$ be an annulus and $f:A'\to H$ be a singular map such that $f(A')=A$ is properly immersed. If $A$ is not parallel into  $\bndry H -\tc$ and the curves $\bndry A$ are essential in $\bndry H$, then there is a properly embedded essential annulus in  $\bndry H-\tc$.
\end{lemma}

\begin{remark}
The proof for this lemma uses a simplified version of the covering space argument used by Freedman, Hass and Scott \cite{fhs2}. The process is more elementary since we are operating in a handlebody.
\end{remark}

\begin{proof}
The basic steps in this proof are to first find another $f$ so that all the lifts of $A=f(A')$ in the universal cover are embedded. We then use subgroup separability to produce a finite sheeted cover of $H$ which contains a lift of $A$ that is an embedded annulus, which has the same properties and does not intersect any of its translates. From this cover we find a regular cover, in which all the lifts of $A$ are embedded. This then implies that the finite regular cover has a non-trivial annulus region and thus so does the original handlebody.

We will assume that the map $f$ is transverse at all times.  Let $G=\pii(H)$, $f_*$ is the induced map on $\pii (A)$ and $f_*(\pii(A))=B\subseteq G$. Therefore $B$ is a free subgroup generated by some $z\in G$.

Let $\bar{H}$ be a cover of $H$ with the projection $\bar{p}:\bar{H}\to H$ such that $\bar p_*(\pii (\bar{H}))=B$. This means there is a lift, $\bar{A}$ of $A$, which is an annulus such that $\pii(\bar H)\cong \pii(\bar A)$. Let $\bar{\tc}=\bar p^{-1}(\tc)$.  As $A$ is not parallel into  $\bndry H - \tc$ then $\bar{A}$ is not homotopic into  $\bndry \br{H} - \bar{\tc}$.

We now want to find an embedded annulus in $\bar{H}$ which is $\pii$-injective and not properly homotopic into $\bndry \bar H - \bar{\tc}$.  Let $n(\bar A)$ be a regular neighbourhood of $\bar{A}$ such that $n(\bar A) \cap \bar{\tc} = \emptyset$.  Then $\br{\bndry \bar{N}\cap \mathrm{int}(\bar{H}) }$ is a set of embedded surfaces.  As $\pii(\bar A) \cong \pii(\bar H)$, we know that there are two of  these  embedded surfaces in $\bar{H}$ whose boundary curves are  essential in $\bar H$. Let one of these be $\bar{A}'$.  Note that as the boundary curves of $\bar A$ are not homotopic  in $\bndry \bar H - \bar \tc$, that is, $\bar A$ is not homotopic into $\bndry \bar H - \bar \tc$, then the simple boundary curves of $\bar A'$ are not homotopic in $\bndry \bar H - \bar \tc$.

By Dehn's lemma and the loop theorem we know that any handles in $\bar{A}'$ can be compressed until $\bar A'$ is an annulus.  After surgering all its handles, $\bar{A}'$ is an essential embedded annulus in $\bar{ \bndry H } -  \bar{\tc}$. Now  let $A=\bar{p}(\bar A')$. As $A$ is compact, it has a finite number of curves of self intersection. We also assume that the self intersection curves are transverse. Let $\bar A_i$'s, for $1\leq i \leq n$, be the lifts of $A'$ in $\bar H$ that intersect $\bar A'$ and $\bar \alpha_i=\bar A'\cap\bar A_i$.  Thus each $\bar \alpha_i$ is a set of singular curves in $\bar A'$.

Let $\tld{H}$ be the universal cover of $H$ and therefore also the universal cover of $\bar{H}$ with the projections $p:\tld H \to H$ and $\tld p:\tld H \to \bar H$, such that $p=\bar p \tld p$.   As $H$ is a handlebody, $\tld H$ is a missing boundary ball,  that is, a ball with a compact set removed from its boundary. As $A$ is $\pii$-injective in $H$, each pullback to $\tld H$ is a  universal cover of $A$, a missing boundary disk.  As $\bar A'$ is embedded in $\bar H$,  each pullback to $\tld H$ is embedded. Then by applying the covering transformation group to $\tld H$ we know that all the lifts of $A$ in $\tld H$ are embedded.

Let $\tld{A}$ be a lift of $\bar A'$ in $\tld H$. Then any lift of $A'$ in $\tld H$, that intersects $\tld A$ must be a lift of one of the $\bar{A}_i$'s in $\bar{H}$. Let $\tld{A}_i$ be some lift of $\bar A_i$ that intersects $\tld{A}$ and $\tld \alpha_i=\tld A \cap \tld A_i$. Note this means that $\tld p (\tld \alpha_i)=\bar \alpha_i$. Also let $\tld{G}$ be the group of deck transformations on $\tld H$ and $\tld{B} \subset \tld{G}$ the stabilizer of $\tld{A}$. Therefore $\tld{G}\cong G$ and $\tld{B}$ is the subgroup of translations along $\tld{A}$. Also let $g_i\in \tld{G}$ where $g_i(\tld{A})=\tld{A}_i$. This means that $g_i \notin \tld{B}$ and that $\tld{B}_i= g_i \tld{B}$ is the set of transformations taking $\tld{A}$ to $\tld{A}_j$.  So for all $b \in \tld B$, $\tld p(b (\tld \alpha_i))=\bar \alpha_i$.

By Hall \cite{ha} we know   there is a finite index subgroup $\tld{L}_i\subseteq\tld{G}$ such that $\tld{B} \subseteq \tld{L}_i$ but $g_i\notin \tld{L}_i$. This property is called subgroup separability. For all $b\in\tld{B}$, $bg_i\tld{A}$ is a translate that intersects $\tld{A}$ and $bg_i\notin \tld{L}_i$. This means that for any $l\in \tld L_i$  that $l(\tld A)\neq b (\tld A_i) = b g_i (\tld A) $ for all $b\in \tld B$. In other words none of the deck transformations in $\tld L_i$ map $\tld A$ to the lift of $\bar A_i$ that intersects $\tld A$.  Let $\hat{H}_i= \tld H/ \tld L_i$ be the cover of $H$ with the fundamental group corresponding to $\tld L_i$ such that $\hat p_i:\tld H \to \hat H_i$.  Therefore $\hat p_i(\tld A)$ is an embedded annulus in $\hat H_i$. Also, for any $b\in B$, $\hat p_i(b \tld A_i) \cap \hat p_i (\tld A) = \emptyset$ and as $\tld L_i$ has finite index in $G$, $\tld {H_i}$ is a finite sheeted  cover of $H$

Therefore $\tld{L}=\tld{L}_1\cap\ldots\cap\tld{L}_n$ is a finite index subgroup of $\tld{G}$ such that for  $l\in\tld{L}$, either $\tld A= l(\tld A)$ or $\tld{A}\cap l(\tld{A})=\emptyset$. Let $\tld H/ L =\hat H$ be the finite sheeted cover of $H$ with the projection $\hat p:\tld H \to \hat H$.  Then $\hat p (\tld A) = \hat A$ is an embedded annulus in $\hat H$ that does not intersect any other lifts of $A'$.

As $L$ has finite index, it must have a finite number of  right cosets, $\{Lx_1, \ldots , Lx_n\}$, for $x_1,\ldots, x_n\in G$. Assume that $Lx_1= L$. Thus if $S_n$ is the group of permutations of $n$ elements, there is a map $\phi : G \to S_{n}$, where $\phi(g)$, for $g\in G$, is the element of $S_n$ that sends $\{Lx_i\}$ to $\{Lx_ig\}$. As both $\phi(g_1)\phi(g_2)$ and $\phi(g_1g_2)$ send $\{Lx_i\}$ to $\{Lx_ig_1g_2\}$, $\phi$ is a homomorphism.  Let $K\subseteq G$ be the kernel of $\phi$. If $g\in K$ then $Lx_i=Lx_ig=Lgx_i$, thus $K\subseteq L$.  As $S_n$ has a finite number of elements, the kernel $K$ is a finite index normal subgroup.  Therefore $\breve H = \tld H/K$ is a finite sheeted normal cover of $H$. Let $\breve p:\breve H \to H$ be the covering projection. Then $\breve H$ is a handlebody and $\breve \tc = \breve p^{-1}(\tc)$ is a set of curves in $\bndry \breve H$ that meet the $n$ disk-condition in $\breve H$.  $\breve H$ is also a cover of $\hat H$ thus all the lifts of $A$ are properly embedded essential annuli in $\bndry \breve H -\breve \tc$.

Then by Freedman Hass and Scott \cite{fhs2}  if we put a Riemannian metric on $H$ and homotop $f$ so that $f(A)$ is least area, then all trivial self  intersections between lifts of $A$ will be removed and thus by  lemma \ref{lemma: intersecting annuli} and lemma \ref{lemma: embedded annuli int vert or horz} all the lifts of $A$ in $\breve H$ are either pairwise disjoint, or intersect each other vertically or horizontally.  If the lifts of $A$ are pairwise disjoint, $A$ must be a properly embedded essential annulus  in $\bndry H -\tc$.  Otherwise  by lemma \ref{lemma:embedded vert int annuli into I-bundle} and lemma \ref{lemma:embedded horz int annuli into tree region} we know that  $\breve H$ must have a non-trivial region $N_I \cup N_T$.  By lemma \ref{lemma: I-bundle in regular cover} and lemma \ref{lemma: tree region in regular cover} we know that the $N_I \cup N_T$ can be isotoped so that its boundary annuli are preserved under $K$ and thus project  to properly embedded essential annuli in $\bndry H -\tc$.
\end{proof}


\begin{lemma}
If $H$ is a handlebody, $\tc$ is a set of triple curves in its boundary that meets the $n$ disk-condition and  $A$ is a properly embedded annulus in $H$ then $A$ is isotopic into $N$.
\end{lemma}

\begin{proof}
Let $A$ be an annulus embedded  in $H$ that cannot be isotoped into $N$.  By lemmas  \ref{lemma:embedded horz int annuli into tree region} and \ref{lemma:embedded vert int annuli into I-bundle} we know that if $A$ has a non-trivial intersection with another embedded annulus  then it can be isotoped into $N_I$ or $N_T$.  Therefore we can isotope $A$ so that it is disjoint from all the boundary annuli of $N$.  This contradicts maximality of  $N$, thus we must be able to isotope $A$ into $N$.
\end{proof}

\begin{lemma}\label{lemma: homotopy of annuli into N}
Let $H$ be a handlebody, $\tc$ be a set of triple curves in its boundary that meets the $n$ disk-condition and $N$ be the annulus region in $H$. If $A$ is an annulus and $f:A\to H - \tc$ is a proper singular essential map then there is a homotopy of $f$  so that $f(A)$ is in $N$.
\end{lemma}

\begin{proof}
To save on notation we will refer to the image of $f(A)$ as $A$ as well.  Let $B$ be the set of boundary annuli of $N$ and $\tc' = \tc \cup \bndry B$.  Then $H'=\overline{H-N}$ is a set of handlebodies such that for any component $H'_j$, the set of essential simple closed curves $\tc'\cap H'_j$  meets the 4 disk-condition in $H'_j$. Also there is a homotopy of $f$ so that $f^{-1}(N)$ is either a set of 4-gons (case 1) or essential embedded annuli (case 2).

Case 1:  All the components of $N$ that $A$ intersects are either in $N_I$ or $N_A$. Assume the singular 4-gons  $H'\cap A$ are essential in $H'$. Then by  the loop theorem we know that there is an embedded essential 4-gon with two boundary arcs in the boundary annuli of $N$.  This contradicts maximality of $N$.

Case 2:  As in the previous case, all the components of $N$ that $A$ intersects are either in $N_T$ or $N_A$. Then by lemma  \ref{lemma: singular annulus} we know that $H'$ must contain an essential properly embedded annulus, contradicting maximality of $N$.

Thus there must be a homotopy of $f$ so that $A$ is disjoint from $B$. If $A$ is not contained in $N$ then  once again by lemma \ref{lemma: singular annulus}, $H'$ contains essential embedded annuli, contradicting maximality of $N$.
\end{proof}


\subsection{Torus theorem}

Let $M$ be a 3-manifold that meets the $(n_1,n_2,n_3)$ disk-condition.  That is for $1\leq i \leq 3$ $H_i\subset M$ is an  embedded handlebody such that; $\bigcup H_i = M$, $\bigcup \bndry H_i = X$ is a 2-complex that cuts $M$ up into the $H_i$'s and $\bigcap H_i =\tc$ is a set of essential simple closed curves that meet the $n_i$ disk-condition in  $H_i$.  We will assume that $(n_1,n_2,n_3)$ is either $(6,6,6)$, $(4,6,12)$ or $(4,8,8)$,  for if the gluing of the three handlebodies meets some disk-condition, it meets one of these three.

\begin{lemma} \label{lemma: homotopy of essential annuli into n-gons or annuli}
Let $M$ be a compact closed 3-manifold that meets the disk-condition as described above.   Suppose $T$ is a torus and $f: T \to M$ is a singular essential map. Then there is a homotopy of $f$ so that either $f(T)$ is disjoint from $\tc$ and for each $i$, $H_i \cap f(T)$ is  a set of essential annuli  or, for each $i$, $H_i \cap f(T)$ is a set of singular disks with essential boundary that intersect $\tc$ exactly $n_i$ times.
\end{lemma}

\begin{proof}
Assume that $f$ is transverse to $X$. Thus $\Gamma = f^{-1}(X)$ is a set of simple closed curves and trivalent embedded graphs in $T$.  Once again let an $(m,n)$-gon be a face of $T$ that is a disk, has $m$ vertices in its boundary and is mapped by $f$ into the handlebody in which $\tc$ meets the $n$ disk-condition.  Let $\Gamma_j$'s be the components of $\Gamma$.  $\Gamma_i$ is a non-essential component if there is a disk $D\subset T$ such that $\Gamma_i\subset  D$.  Then by lemma \ref{lemma: removing pullback graphs} we know that there is a homotopy of $f$ to remove $\Gamma_i$ and hence remove all non-essential components of $\Gamma$.

Therefore there are two cases.  Either all faces of $\Gamma$ are disks or $\Gamma$ has faces which are essential annuli.  Note that $f(T) \cap X \not= \emptyset$ as $f$ is $\pii$-injective and $\pii (H_i)$ doesn't have a free abelian subgroup of rank 2.

If $\Gamma$ is connected then all the faces must be $(m,n)$-gons and all the vertices have order three.  Let $\mathbf F$ be the set of faces of $T$.  We can then put a metric on $T$, as we did in the proof of lemma \ref{lemma: removing pullback graphs}.  So all the edges are geodesics of unit length and if $F\in \mathbf F$ is an $(m,n)$-gon then the angle at each vertex is $\pi(1-2/n)$ and there is a cone point in $int(F)$.  Once again this means that the curvature around each vertex is $2\pi$.  Let $K(F)$ be the curvature at the cone point in $F$. By the Gauss-Bonnet theorem we know that $$\mathbf K(F)  = 2\pi(1-m/n).$$

Therefore if $m>n$ then $\mathbf K (F) < 0$, if $m=n$, $\mathbf K (F) = 0$ and if $m<n$, $\mathbf K (F) > 0$.  Also by the Gauss-Bonnet theorem we know that  $$\sum_{\mathbf F} \mathbf K(F) = 0$$

Therefore if $\mathbf F$ contains an $(m,n)$-gon such that $m>n$, then it must also contain a face $F$ such that $m<n$. Thus by the disk-condition we know that $f (\bndry F)$ is not essential in $\bndry H_k$.  So $f(F)$ is parallel into $\bndry H_k$ and there is a homotopy of $f$  so that $f(F)\subset \bndry H_k$.  We can then push $F$ off $\bndry H_k$ removing the face $F$ from $\mathbf F$. Note that when we do this, the order of the faces adjacent to $F$ either decreases by two or an $(m,n)$-gon and an $(m',n)$-gon merge to become an $(m+m'-4,n)$-gon, as shown in figure \ref{fig:homotop f to remove gon.}. We can repeat this process as long as $\mathbf F$ contains faces with positive curvature.  Each time we do this move, we reduce the number of faces in $\mathbf F$ by at least one.  Therefore this process must terminate after a finite number of moves, when all the faces are $(m,n)$-gons such that $m=n$.

Now let's look at the case where $\Gamma$ contains more than one component.  Let $\Gamma_i$ be a component of $\Gamma$.  Then $\Gamma_i$  cuts $T$ up into faces that are a single annulus and a sum of a number of disks.  Let $A$ be the union of $\Gamma_i$ and the faces which are disks.  Now we know that the Euler characteristic of $A$ is 0.  Put a metric on $A$ as we did above.  $\Gamma_i$ must have boundary vertices,  that is vertices adjacent to less than three faces of $A$.   Thus using the same arguments using the Gauss Bonnet theorem we know that $A$ must have some face with positive curvature.  This means that such faces are boundary parallel in the handlebody and there is a homotopy of $f$ to remove them.  As before this process can be repeated until all the components are simple closed essential loops.
\end{proof}

We are now ready to prove the torus theorem.

\emph{proof of theorem \ref{thrm: torus}.}
Let $N_i$ be the maximal annulus region for $H_i$ and $P_i$ be the maximal disk region for $H_i$.  The idea of this proof is to show that we can find  submanifolds of either the $N_i$'s or the $P_i$'s such that when we glue them together, the resulting embedded sub-manifold can be fibered by $\s[1]$ and either has essential tori boundary or the fibering can be extended to the whole of $M$.  In the interest of reducing notation, the image of $f(T)$ in $M$ will be denoted as $T$.  Thus when we talk about a homotopy of $T$, we are implying a homotopy of $f$.

By lemma \ref{lemma: homotopy of essential annuli into n-gons or annuli} there is a  homotopy so that  either $T$ is disjoint from $\tc$ and for each $i$, $H_i \cap T$ is  a set of essential singular annuli  not homotopic into $\bndry H - \tc$  or, for each $i$, $H_i \cap T$ is a set of singular meridian  disks that intersect $\tc$ exactly $n_i$ times.

The first case is when $T$ is disjoint from the triple curves and $H_i \cap T$ is a set of singular essential annuli.  We can also assume that no components of  $H_i\cap T$ are parallel into $\bndry H_i - \tc$.   By lemma \ref{lemma: homotopy of annuli into N} we can isotope each $N_i$ so that  $H_i\cap T\subset N_i$.

Let $A_i = X \cap N_i$, then $A_i$ is a set of essential surfaces in $\bndry H_i$. Note that  $T\cap \bndry H_i \subset A_i$ and thus $T\cap X \subset \bigcup_{i\not= j} (A_i\cap A_j)$. We will first reduce $N_1$.  Let $S_i = A_i \cap (A_j\cap A_k)$, where $i$, $j$ and $k$ are different.  Let $N_1'$ be the maximal subset of $N_1$ such that $N_1'\cap X \subseteq S_1$ and each component of $\overline{\bndry N_1' \cap int(H_1)}$ is an essential annulus. There are three cases corresponding to components of $N_I$, $N_T$ and $N_A$.

Let $B$ be a component of $N_1$ such that $B$ is an $I$-bundle region and $F$ is its base space. Then let  $F'\subseteq F$ be the maximal sub-surface such that $B'\cap \bndry H_1 \subseteq S_1$, where $B'$ is the $I$-bundle over $F'$. Then $B'$ is a component of $N'_1$. Note that components that do not intersect $S_1$ are removed.

If $B$ is a tree region then it is a fibered solid torus and $B\cap \bndry H_1$ is a set of essential annuli. Then there is an isotopy of $B$ such that each annulus in $B\cap \bndry H_1$ is either contained in $S_1$ or  in $int(H_1)$.  Note that some annuli in $\bndry H_1$ may get pushed into $int(H_1)$. Let $B'$ be the resulting fibered torus. Note that when the number of annuli in $B\cap \bndry H$ is reduced to produce $B'$, the fibering of the torus is still parallel to the boundary curves of the boundary annuli.  Then $B'$ is a component of $N'_1$.  If $B'\cap H_1 = \emptyset$ we remove it from $N'_1$.

If $B$ is a component of $N_A$, as defined in section \ref{section: Annulus region} then either it can be isotoped so that $B\cap H_1\subseteq S_1$ or it is removed.  As  $T\cap X \subset \bigcup_{i\not= j} (A_i\cap A_j)$ we know that $N_1'\not= \emptyset$.  We now let $N_1=N_1'$.

We now repeat this process for each $N_i$ in turn until the process stabilises. That is for $i\not= j$, $i \not= k$ and $k \not= j$, $A_i = \bndry H_i \cap (A_j \cup A_k)$.  We know that it stabilises before $\bigcup N_i = \emptyset$ because $T\subset  \bigcup N_i$.

Next we want to change the fiberings of the $N_i$'s so that all components that are regular neighbourhoods of embedded annuli or Mobius bands are fibered by $\s[1]$. This means that for any component $B$ of $N_i$ such that $B\cap \bndry H_i$ is a set of annuli, then $B$ is a fibered solid torus, otherwise it is an $I$-bundle.  Now when we let $\mathbf N= \bigcup N_i$ and  all the fiberings of components match, then $N$ is a Seifert fibered sub-manifold of $M$ and  $\bndry \mathbf N$ is a set of embedded tori.

By lemma \ref{lemma: compressing disk of i-bundle has int at least n/2}, if $\mathbf N_j$ is a  component of $\mathbf N$ such that $H_i \cap \mathbf N_j$ is an $I$-bundle with an base space that is not an annulus or a Mobius band, then the boundary tori of $\mathbf N_j$ are essential in $M$. The final step in this case is to either make all the boundary tori of $\mathbf N$ essential or expand $\mathbf N$ so that $\mathbf N = M$.  If $\mathbf N_j$ is a component of $\mathbf N$ and $F \subset M $ is an embedded solid torus such that $\bndry F \subseteq \mathbf N_j$,  then either $F\cap \mathbf N_j = \bndry F$ or $F\cap \mathbf N_j = \mathbf N_j$.  If $F\cap \mathbf N_j = \bndry F$ we then add $F$ to $\mathbf N$ and extend the fibering to it. This can always be done as the fibers of the component are essential in $M$. Therefore the meridian disk of the solid torus being added cannot be parallel to the fibering of $\mathbf N_j$.   If $\mathbf N_j$ is contained in $F$ we remove $\mathbf N_j$ from $\mathbf N$.  This  process is repeated until either all boundary tori are essential or $\mathbf N = M$.  We know the process will terminate before all of $\mathbf N$ has been removed because $T\subset \mathbf N$ and $T$ is essential. Thus the component containing $T$ cannot be contained in a solid torus.

The next case we look at is where $H_i\cap T$ is a set of singular $n_i$-gons. Let $ P_i$ be the disk region in the handlebody $H_i$. Next we want to define a process for reducing components of  $ P_i$ until all their boundaries coincide in $X$ and then show that we can expand the `core' fibering to the whole sub-manifold.   Let $A_i = X\cap  P_i$.  By lemma \ref{lemma: singular disks into C} we know that we can isotope each $ P_i$ so that $H_i\cap g(T)\subset  P_i$.  Thus $T \cap \bndry H_i \subset P_i\cap(P_j \cup P_k)$, for $i \not= j$, $j\not= k$ and $k\not= i$.

Reduce $P_1$ so that $P_1\subseteq P_2\cup P_3$.  By reducing we mean chop off fingers that don't match up, reduce base spaces of  the cores and possible remove entire components of $P_1$.  This process finishes before $P_1$ is entirely removed as   $T \cap \bndry H_i \subset P_i\cap(P_j \cup P_k)$.  Note that if a component of $P_1$ is reduced to the regular neighbourhood of a single meridian disk we  forget the fibering of its core.  As we reduce $P_1$, $\overline{ \bndry P_1 \cap int(H)}$  remains a set of  essential annuli and meridian disks.

This process is repeated in turn for each $P_i$.  Once again we know that the process stabilises before all the $P_i$'s  are removed as $T \cap \bndry H_i \subset P_i\cap(P_j \cup P_k)$.  All the components with fibered cores obviously match up to be  fibered tori in $\mathbf P = \bigcup P_i$.  Clearly these do not intersect so the fibering can be extended across $\mathbf P$. $\mathbf P$ is a Seifert fibered submanifold of $M$ and  each of the boundary tori of $\mathbf P$ is tiled by either meridian disks or essential annuli that are essential in $T$.  As before if any of the torus boundaries of $\mathbf P$ are not essential, they are either filled in with a solid torus or removed.
\qed


\subsection{Characteristic variety}

Finally we wish to show that both flavours of characteristic variety fit together nicely. That is, if the flavours intersect, their $\s[1]$ fiberings can always be made to agree. If either component is a $T^2 \cross I$ this is easy.  Thus we want to study the case where each component has a unique fibering.

Let $\mathbf N$ be the maximal annulus region in $M$ and $\mathbf P$ be the maximal disk region.  By the usual arguments we can see that both are unique up to isotopy. We can also assume that $\mathbf N$ is disjoint from $\tc$ and that both flavours have non-empty boundary. Thus $\bndry \mathbf N \cup \bndry \mathbf P$ is a set of essential embedded tori. If $\mathbf N\cap \mathbf P = \emptyset$, then there is no problem. Therefore we can assume that $\mathbf N \cap \mathbf P \not= \emptyset$.  Let $N'$ be a component of $\mathbf N$ and $P'$ be a component of $\mathbf P$ such that $N' \cap P' \not= \emptyset$.  It is not possible for $P' \subset N'$ and if $N'\subset P'$ there is no problem.  Therefore we can assume that there is a boundary torus $B\subset \bndry P'$ such that $B\cap N' \not= \emptyset$.  As $\bndry N'$ is a set of essential tori, $B\cap N' $ is a set of essential annuli in $N'$. Thus $H_i\cap(B\cap N')$, for any $i$, is a set of quadrilaterals.  Therefore, if the components of $H_i \cap N'$ are fibered by $\s[1]$, then $N' \homeo T^2 \cross I$. Thus we can assume that $N'$ is fibered such that $N'\cap H_i$ is a set of $I$-bundles.  Therefore it just remains to show that $H_i\cap(N'\cap P')$ is an $I$-bundle.

Let $F$ and $F'$ be two meridian disks in $H_i$ that intersect $\tc$ $n_i$ times and have a non-trivial intersection and  $A$ be an essential properly embedded annulus in $H_i - \tc$. We can  assume that $A$ has been isotoped so that  $F \cap A$ is a set of disjoint properly embedded  arcs in $F$. If any of the arcs in $F\cap A$ are not bisecting then $A$ is boundary parallel.  In this case $F'\cap A$ cannot contain any properly embedded arcs, for if it did, this would provide an isotopy of $F$ to remove that intersection between  $F$ and $F'$.  Thus $F\cap A$ must be a set of bisecting arcs in $F$, similarly $F'\cap A$ is a set of properly embedded bisecting arcs in $F'$ and $A$ is not boundary parallel.  If we then let $Q$ be the regular neighbourhood of $F \cup F'$, then $B = \overline{\bndry Q \cup int(H)}$ is a set of properly embedded annuli and meridian disks that intersect $\tc$ exactly $n_i$ times.  As in the proof of lemma  \ref{lemma: no triple points in intersecting n-gons}, there is an isotopy of $A$ so that $A\cap B$ is a set of properly embedded parallel arcs that are not boundary parallel in $A$.  Thus there is an isotopy to remove any triple points.

The components of $P'\cap H_i$ can be thought of as regular neighbourhoods of a set of meridian disks that intersect $\tc$ exactly $n_i$ times. From above, if there are two meridian disks in $H_i$ that intersect $\tc$ $n_i$ times and have a non-trivial intersection, then any essential annulus can be isotoped so that it is disjoint from their intersection.  Lemma \ref{lemma: compressing disk of i-bundle has int at least n/2} says any boundary compressing disk of the annuli $N' \cap H_i$ has order at least $n_i/2$. Therefore the intersection between boundary annuli of $N' \cap H_i$  and a meridian disk of order $n_i$ must be bisecting in the meridian disk.  By these two observations we can see that $H_i\cap (N'\cap P')$ is an $I$-bundle.

\subsection{Atoroidal manifolds}

An interesting question raised by Cameron Gordon is to find an additional condition that would result in this class of manifolds being atoroidal. By lemma \ref{lemma: homotopy of essential annuli into n-gons or annuli}, a sufficient condition for a manifold $M$ that meets the disk-condition to not contain any essential tori that intersect the triple curves, is the manifold meets a stronger disk-condition with $\sum 1/n_i<1/2$.  A sufficient condition that $M$ does not contain any essential tori disjoint from the triple curves is that in at least two of the handlebodies, any essential annuli disjoint from $\tc$ are boundary parallel.

Let $H$ be a handlebody and $\tc$ an essential set of disjoint simple closed curves in $\bndry H$ that meet the $n$ disk-condition.  Let $A$ be a properly embedded essential annulus in $H$ disjoint from $\tc$. Then by lemma \ref{lemma: waveless minimal system of m-disks}, $H$ has a waveless minimal system of disks, $\mathbb D$, see definition \ref{defn: waveless system}.  Let $B$ be the 3-ball produced when $H$ is cut along $\mathbb D$, $S\subset \bndry B$ be the punctured sphere produced when $\bndry H$ is cut along $\mathbb D$ and $\Gamma = \tc \cap S$.  As in the proof for lemma \ref{lemma: waveless minimal basis of m-disks}, let $\Gamma'\subset \s[2]$ be the graph produced by letting components of $\bndry S$ correspond to vertices and parallel components of $\Gamma$ correspond to single edges, see figure \ref{fig: triple curve to graph.}.

As $A$ is a properly embedded essential annulus, $B\cap A = \{A_1, ... , A_k\}$ is a set of properly embedded quadrilaterals in $B$ such that, for any $i$, $A_i \cap S$  is two properly embedded arcs in $S$. An equivalent statement to $A$ being boundary parallel is that the curves $\bndry A$ are parallel in $\bndry H$ or that for each $i$, the arcs $A_i \cap S$ are parallel in $S$.

\begin{lemma}\label{lemma: annuli boundary parallel}
    If $\Gamma'$ is maximal and contains no 2-cycles (definition  \ref{defn: 2-cycle}) then all properly embedded annuli in $H$ disjoint from $\tc$ are boundary parallel.
\end{lemma}

\begin{proof}
By maximality of $\Gamma'$, the arcs of $A_i\cap S$, for all $i$, must be parallel to some arc of $\Gamma$ and  as $\Gamma'$ contains no 2-cycles, both arcs of $A_i \cap S$ must be parallel to the same arc of $\Gamma$ and thus parallel. Therefore, from above, any properly embedded essential annulus in $H - \tc$ must be boundary parallel.
\end{proof}

Let $K\subset \s[3]$ be  an $(a_1,a_2,a_3)$ pretzel link such that, for each $i$,  $a_i \geq 4$ and the spanning surface $F$ shown in figure \ref{fig:(333)pretzel_knot.} is orientable.  As in section \ref{subsection: dehn filling example}, let $M$ be the manifold produced by taking the 3-fold branched cover of $\s[3]$ with $K$ as the branch set and $X$ be the 2-complex produced by gluing the lifts $F$ in $M$. Then $M$ meets the disk-condition and $X$ is a 2-complex that cuts it up into injective handlebodies.  As $a_i\geq 4$, the basis bounded by the curves shown in figure \ref{fig:(333)pretzel_knot.} is an 8-waveless basis (definition \ref{defn: n-waveless}) for $K$ in the handlebody $\overline{\s[3] - S}$.  Therefore  all meridian disks in the handlebody $\overline{\s[3] - S}$ intersect $K$ at least 8 times. We can produce a waveless minimal system of meridian disks for the handlebody $\overline{\s[3] - f}$ by removing any one of the disks from the basis. The associated graph $\Gamma'$, as constructed above meets the conditions of lemma \ref{lemma: annuli boundary parallel}.  Thus the 3-fold branched cover of such a pretzel link meets the disk-condition and is atoroidal.

\bibliographystyle{plain}
\bibliography{bibfile}

\begin{thebibliography}{10}

\bibitem{Ai&Ru1}
I.~R. Aitchison and J.~Hyam Rubinstein.
\newblock Localising {D}ehn's lemma and the loop theorem in 3-manifolds.
\newblock {\em Math. Proc. Cambridge Philos. Soc.}, 137(2):281--292, 2004.

\bibitem{fhs1}
M.~Freedman, J.~Hass, and P.~Scott.
\newblock Closed geodesics on surfaces.
\newblock {\em Bull. London Math. Soc.}, 14:385--391, 1982.

\bibitem{fhs2}
M.~Freedman, J.~Hass, and P.~Scott.
\newblock Least area incompressible surfaces in {$3$}-manifolds.
\newblock {\em Invent. Math.}, 71(3):609--642, 1983.

\bibitem{ha}
M.~Hall.
\newblock Coset representation in free groups.
\newblock {\em Trans. Amer. Math. Soc.}, 67:421--432, 1949.

\bibitem{Hem1}
J.~Hempel.
\newblock {\em {$3$}-{M}anifolds}.
\newblock Princeton University Press, Princeton, N. J., 1976.
\newblock Ann. of Math. Studies, No. 86.

\bibitem{ja1}
W.~Jaco.
\newblock {\em Lectures on three-manifold topology}, volume~43 of {\em CBMS
  Regional Conference Series in Mathematics}.
\newblock American Mathematical Society, Providence, R.I., 1980.

\bibitem{Mat1}
S.~Matveev.
\newblock {\em Algorithmic topology and classification of 3-manifolds},
  volume~9 of {\em Algorithms and Computation in Mathematics}.
\newblock Springer-Verlag, Berlin, 2003.

\bibitem{Oe1}
U.~Oertel.
\newblock Closed incompressible surfaces in complements of star links.
\newblock {\em Pacific J. Math.}, 111(1):209--230, 1984.

\bibitem{ro1}
D.~Rolfsen.
\newblock {\em Knots and links}.
\newblock Publish or Perish Inc., Berkeley, Calif., 1976.
\newblock Mathematics Lecture Series, No. 7.

\bibitem{rou&sa1}
C.~P. Rourke and B.~J. Sanderson.
\newblock {\em Introduction to piecewise-linear topology}.
\newblock Springer-Verlag, New York, 1972.
\newblock Ergebnisse der Mathematik und ihrer Grenzgebiete, Band 69.

\bibitem{sc}
P.~Scott.
\newblock There are no fake {S}eifert fibre spaces with infinite {$\pi
  \sb{1}$}.
\newblock {\em Ann. of Math. (2)}, 117(1):35--70, 1983.

\end{thebibliography}
\end{document}